\documentclass[10pt,twocolumn]{article}

\usepackage[utf8]{inputenc}
\usepackage[T1]{fontenc}
\usepackage[margin=0.75in]{geometry}

\usepackage{hyperref}
\usepackage{url}
\usepackage{booktabs}
\usepackage{amsfonts}
\usepackage{nicefrac}
\usepackage{microtype}
\usepackage{graphicx}
\usepackage{doi}
\usepackage[cmex10]{amsmath}
\usepackage{amssymb}
\usepackage{pdflscape}
\usepackage{longtable}
\usepackage{changepage}
\usepackage{cite}
\usepackage{float}
\usepackage[table]{xcolor}
\usepackage{array}
\usepackage{caption}
\usepackage{cleveref}
\usepackage{authblk}

\definecolor{llnlLivermoriumIce}{gray}{0.95}
\definecolor{llnlElementalNavy}{RGB}{0,40,85}
\definecolor{llnlEnergeticAzure}{RGB}{51,102,204}
\definecolor{llnlImpactBlue}{RGB}{0,50,161}

\newcommand{\rothead}[1]{\rotatebox[origin=c]{90}{\parbox{2.2cm}{\centering #1}}}

\title{High-Fidelity Capacity Expansion Planning for Puerto Rico's Electric Power System} 
 
\author[1]{Elizabeth Glista}
\author[1]{Tomas Valencia Zuluaga}
\author[1]{Amelia Musselman}
\author[1]{Minda Monteagudo}
\author[1]{Juliette Franzman}
\author[1]{Jean-Paul Watson}

\affil[1]{Lawrence Livermore National Laboratory}

\date{June 26, 2026}

\hypersetup{
pdftitle={Capacity Expansion Planning for Puerto Rico's Electric Power System},
pdfauthor={Elizabeth Glista, Tomas Valencia Zuluaga, Amelia Musselman, Minda Monteagudo, Juliette Franzman, Jean-Paul Watson},
pdfkeywords={capacity expansion planning, Puerto Rico, power systems, unit commitment, stochastic optimization, generation planning},
}

\newcommand{\myth}[1]{\color{white}\textbf{#1}}

\colorlet{RowOdd}{white}
\colorlet{RowEven}{gray!10}

\begin{document}

\twocolumn[
\maketitle
\begin{abstract}
This study presents a mathematical optimization framework and analysis to inform practical long-term investment planning in Puerto Rico's electric power system. In this study, we utilize a high-resolution capacity expansion planning (CEP) model to identify least-cost generation and storage investments that significantly improve the reliability of Puerto Rico's power system. The model co-optimizes new investments along with thermal generator retirements, while considering a detailed operational model of generator dispatch, thermal generator commitments, fuel selection, and storage operation.  Constraints are imposed to enforce equipment engineering limits, system constraints, fuel supply limitations, and load satisfaction. Key modeling advances of this study, relative to prior long-term planning studies on Puerto Rico's power system, include: (i) a nodal representation of the transmission network at 38 kV and above; (ii) hourly chronological simulation of system operations for a large set of representative days; (iii) an explicit unit commitment formulation for all existing and new thermal units, with realistic ramping, minimum up- and down-times, and startup costs; (iv) explicit modeling of regional fuel supply constraints; and (v) operational scenarios that capture load variability, renewable availability, and the high forced outage rates of legacy units.\\

Using data supplied by LUMA, the Puerto Rico Electric Power Authority (PREPA), and the U.S. Department of Energy (DOE), along with publicly available data, we construct present day (2024) and future (2030) Puerto Rico representative systems, where the future system includes planned generation and storage projects. We then analyze several planning scenarios with varying future load levels, fuel supply assumptions, materializations of currently planned additions to the generation and storage fleets, and allowable new technologies. Across the scenarios examined and under study assumptions, we found that least-cost portfolios that maintain modeled reliability targets include approximately 1.5 GW or more of new H-class combined cycle (CC) capacity in addition to currently planned projects.
These additions are needed primarily to replace unreliable legacy thermal units rather than to serve new load. The new CC investments eliminate modeled load shedding in the bulk system and restore a robust reserve margin, even under stressed load and outage conditions.\\

\end{abstract}


\vspace{1em}
]

This work was performed under the auspices of the U.S. Department of Energy by Lawrence Livermore National Laboratory under Contract DE-AC52-07NA27344. The work was supported by the Department of Energy (DOE) Office of Electricity's NAERM program and the LLNL-LDRD Program under Project No. LDRD25-SI-007 (LLNL-JRNL-2017618). The views expressed in the article do not necessarily represent the views of the DOE or the U.S. Government.



\section{Introduction}

Puerto Rico's power grid is an isolated system serving approximately 1.5 million customers \cite{puerto_rico_electric_power_authority_puerto_nodate}. The territory's island geography (in an area prone to hurricanes) coupled with decades of operational and financial challenges have led to a rapidly deteriorating situation over recent years.
The frequency and duration of outages in Puerto Rico have increased since 2021, even in the absence of major weather events, with customers experiencing on average up to 13.7 interruptions annually, totaling approximately 30 hours of interruptions per year \cite{aramayo_even_2025}, over twice the national average and far above other island systems like Hawaii or hurricane-prone areas like the Florida Keys \cite{enernex_enernex_nodate}.
Among the challenges contributing most to Puerto Rico's reliability deficiencies are an outage-prone transmission grid and an aging generation fleet with frequent planned and unplanned outages. As a result, operating with insufficient system reserves has become frequent in Puerto Rico \cite{luma_monthly_2025}.

Several efforts have been undertaken by local and federal authorities to bring the quality of the electricity service in Puerto Rico up to U.S. national standards, while simultaneously pursuing other objectives like the goal, enacted into law by the Puerto Rico legislature in 2019, of meeting 100\% of its electricity needs with renewable sources by 2050.
Some of the actions taken during this time included the arrival of LUMA, a new entity tasked with operating the power system and coordinating the Integrated Resource Planning (IRP) process \cite{luma_integrated_2024}. Under this framework, a number of storage and generation projects, mostly solar photovoltaic and wind power plants, have been announced and find themselves in various stages of development. These projects are organized into so-called ``tranches,'' a term we sometimes use in this text.
After about five years with limited success and a situation that has continued to deteriorate, recent measures taken include a short-term (12-36 months) stabilization plan \cite{nepr_stab_plan_2025} to identify projects that must be prioritized, notably: postponing the retirement of existing coal-fired generation, extending and adding short-term emergency generation units, and confirming and adding battery storage capacity to the bulk power system.
Recent developments also include relaxing the interim renewable target goals in the interest of improving reliability for  customers \cite{negociado_de_energia_de_puerto_rico_competitive_2025} and moving to change the current system operator \cite{acosta_vilanova_i_2025}.

The U.S. Department of Energy (DOE) has led various efforts to mobilize the U.S. National Laboratories and provide modeling support to Puerto Rico.
Most relevant for the work presented in this report was the PR100 study \cite{baggu_puerto_2024}, in which different pathways are investigated to reach the 100\% renewable target by 2050, while improving reliability, resilience, and system performance. 
This extensive study spanned data gathering and generation, modeling, and a significant social component. The modeling included production cost modeling, capacity expansion planning, resource adequacy studies, and power flow validation. 

The subsequent work performed by Lawrence Livermore National Laboratory (LLNL) and presented here is much less extensive in scope than the PR100 study and focuses only on capacity expansion planning (CEP).  
Our CEP model differs from the previously mentioned study's model in two important aspects.  First, we do not treat the 100\% renewable target as a requirement but consider a renewables-only scenario as only one possible planning scenario among many in the study. Second, we consider more diverse and higher-fidelity thermal unit representations, specifically natural gas units. 
On the methodological front, a few key aspects differentiate this study from previous modeling efforts for Puerto Rico.
Most importantly, having identified that the operational limitations and frequent outages of the legacy thermal units in the system have a substantial impact on system reliability, we expended significant effort into accurately capturing the operational details of these units, which are typically ignored in long-term planning studies.
    This includes detailed modeling of the unit commitment problem as part of the CEP model. Other differences include:
\begin{itemize}
\item Natural gas fuel supply constraints are explicitly modeled.
\item Stochastic operational scenarios including generator outages are included in the CEP model.
\item The temporal resolution is hourly across all the operational scenarios considered.
\item The system is modeled with nodal resolution at voltage levels of 38kV and above, as opposed to only 230kV and above.
\item Generation units are forced to be exact multiples of real turbine sizes, as opposed to allowing fractional values during optimization (and subsequently rounding to restore integrality) -- which may underestimate investment costs.    
\end{itemize}

We note that the PR100 study included a capacity expansion model with phased decisions (2025-2050), considering the path dependency of those decisions, and included more renewable generation technologies (including longer duration storage and biomass) than we consider in this paper.  Naturally, these approaches may produce different results.

\subsection{Scope of the Study}

This study serves as an exploratory analysis to determine a least-cost investment plan in generation and storage infrastructure for the island of Puerto Rico based on a set of exogenous scenarios characterizing fuel availability, technology portfolios, and target operational years.
The selected scenarios are intended to illustrate how investment recommendations and costs change under different assumptions about load, fuel supply, planned projects, and allowable technologies.  They are not precise forecasts or expressions of policy preferences.

The high-fidelity modeling of the operational limitations of thermal power plants and the high cost attributed to the involuntary disconnection of customers (i.e., load shedding) are intended to enable identification of least-cost investment portfolios that improve modeled service reliability, subject to the operational and infrastructure constraints imposed in the study.
Setting the objective for our modeling tool as total-cost minimization with a high penalty for unserved energy in this way is expected to result in an investment plan that is in the best interest of the electricity users, both in terms of reliability of service and service cost.
Nevertheless, a resource adequacy (RA) study 
and a detailed analysis of the implications of any proposed investment plan for the  ratepayer's electricity bill were out of the scope of this study. They should be performed as part of the modeling work used to quantify decisions made in the Puerto Rico infrastructure rebuild effort.

The CEP model used here makes significant simplifications regarding the constraints imposed by the real transmission system. To partially address this limitation, a selection of the investment plans obtained from our model was analyzed by a team at Pacific Northwest National Laboratory (PNNL). They performed power flow, contingency, and dynamic stability analyses on the Puerto Rico system that would result from the implementation of these investment plans. The feasibility of the proposed system was assessed, and in particular, the need for additional generation and/or battery storage reserves for frequency and voltage support.
While preliminary results show that the investment plans obtained from the CEP model presented in this study are operationally feasible and reliable when analyzed with higher-fidelity operational models, exhaustive analysis in this direction was out of the scope of this study.

The remainder of this report is organized as follows. 
In Section \ref{sec:cep}, we introduce our Capacity Expansion Planning model and describe in detail the mathematical formulation of this problem as an optimization model as well as modeling assumptions and limitations.
In Section \ref{sec:test_system}, we describe the data sources used as inputs for our CEP model of Puerto Rico and assumptions made about missing or incomplete data to enable this analysis.
The different planning scenarios and corresponding results are presented in Section \ref{sec:simulations}. The conclusions of the study are presented in Section \ref{sec:conclusions}.
The Appendix contains details about the data that were relegated from the main part of this text.

\section{Capacity Expansion Planning}
\label{sec:cep}
To determine a cost-optimal long-term investment plan for the Puerto Rico power system, we consider the problem of Capacity Expansion Planning (CEP), a structured mathematical optimization problem that can be formulated as a stochastic mixed-integer linear program (MILP).  The model aims to find the best portfolio of new generators and storage, and the sizes and locations of these units, subject to constraints on land-use, fuel availability,  thermal limits of transmission lines, thermal generator commitments, generator availability, and demand projections.  The objective of the model is to minimize overall system investment and operational costs, in addition to a high penalty on load shedding to encourage high system reliability.  From a planning perspective, the purpose of this formulation is not only to identify least-cost capacity additions but also to test whether candidate portfolios remain operationally feasible under realistic conditions of particular importance in the Puerto Rico system, including hourly variability, unit commitment constraints, fuel limitations, and elevated outage rates for legacy thermal units.

There are three key features of our model of particular relevance: (1) the model is a nodal model with high spatial and temporal resolution; (2) it includes a high-resolution representation of the operational details of thermal generators, including constraints on their fuel supply and unit commitment; and (3) it captures weather variance through numerous representative day samples.  
While CEP is a well-studied problem in the power systems literature, most existing models are zonal models with low spatial resolution, e.g., using 34 buses to model the entire Western Interconnection of the United States \cite{krishnan_evaluating_2016}.  In our recent case study \cite{glista2025}, we found that lower spatial resolution models led to 36\% under-investment in new generation and corresponded to losses of around 40\%, equivalent to billions of dollars a year for a realistic test system.  Recent work has also demonstrated how higher-resolution modeling of thermal generator commitments can impact CEP results, e.g., by encouraging investments in more flexible and dispatchable resources \cite{palmintier_impact_2016,poncelet_unit_2020}.  High temporal resolution allows for higher fidelity modeling of renewable power sources and thermal generator commitments \cite{mallapragada_impact_2018,marcy_comparison_2022}.
Lastly, because of the nature of electric grid planning, there is a high amount of variability associated with some of the model's inputs, in particular, demand for electric power and the availability of non-dispatchable types of generation (e.g., wind and solar).  Thus, it is critical to consider stochastic planning models that can incorporate this variability via a large set of possible scenarios.  There is significant evidence that stochastic models should be used instead of deterministic ones for CEP \cite{Scott2021, munoz_2014}.  In this study, our stochastic model captures variability from load, thermal generator outages, and the availability of wind and solar.  This stochastic framework could also be used to capture longer-term uncertainty rather than just hourly/daily variability, including uncertainty in fuel prices and load growth scenarios.  For the scope and interpretability of this study, these are treated as different sensitivities of the model, rather than included as scenarios in a already very large stochastic problem.

Our CEP model is a high-resolution nodal model that co-optimizes new generation and grid-scale storage investments in a power system.  Although the general model does support transmission expansion, we do not consider transmission expansion here due to limited data availability for Puerto Rico. Our model is stochastic in that we consider several representative scenarios (modeled via representative days) of electric power demand, wind and solar availability, and forced outages of thermal generators.  These representative days capture the variability inherent in these data inputs to the model.
Our model an adaptation of the formulation introduced in \cite{Go2016} and further developed in \cite{zuluaga_parallel_2024, Musselman2024}.  Unlike these models, we also consider constraints from the commitment of thermal generators, including ramping constraints and startup/shutdown costs and constraints, as well as more advanced fuel supply modeling, including dual fuel switching and area-specific fuel supply constraints.  Our formulation includes both investment-stage and operational-stage costs and constraints, where the operational-stage problems are effectively modified Unit Commitment (UC) / Economic Dispatch (ED) problems \cite{knueven2020mixed}.  The model includes constraints that link the second-stage decisions to the investment-stage decisions, e.g., the power generation at a bus must be bounded by the amount of installed generation capacity at that bus. A detailed formulation is provided below.

\subsection{Notations}
The symbols $ \mathbb{R} $ and $ \mathbb{Z} $ denote sets of real continuous and integral numbers, respectively. The symbol $ \mathbb{R}^N $ denotes the space of $N$-dimensional real vectors.  The symbol $\mathbb{R}_+$ denotes the space of real numbers with non-negative values. The symbol $| \cdot |$ is the absolute value operator if the argument is a scalar or vector; otherwise, it is the cardinality of a measurable set.  Transmission branches consist of transmission lines and transformers.

\subsection{Problem Setup}

Given some initial power network defined as a set of buses $\mathcal{B}$, transmission branches $\mathcal{L}$, existing generators $\mathcal{G}^0$, and existing storage elements $\mathcal{S}^0$, we take the set of buses and branches to be fixed, and simultaneously consider the decisions of 1) building new generators, 2) building new grid-scale storage elements, and 3) retiring existing thermal units that are not bound by contractual commitments.  For now, we are ignoring transmission expansion, but this could be considered in future work.  

To formalize the model, we consider the set of candidate thermal generators to be $\mathcal{G}_T^*$, candidate renewable generators to be $\mathcal{G}_R^*$, and candidate storage elements to be $\mathcal{S}^*$.  Note that we classify all candidate generators as either ``thermal'' or ``renewable'' such that the set of all candidate generators is given by the union of the two disjoint sets $\mathcal{G}^*\triangleq \mathcal{G}^*_T\cup \mathcal{G}^*_R$.  This distinction is included for the modeling of thermal generator operation, not for any policy-specific aims, i.e.,  our model is technology-agnostic unless otherwise noted.  We note that the ``renewable'' generator category includes both dispatchable (e.g., hydropower) and non-dispatchable (e.g., wind and solar) power sources.

The set of all generators, both existing and candidate, is given by $\mathcal{G}\triangleq \mathcal{G}^0\cup \mathcal{G}^*$, and the set of all storage elements, both existing and candidate, is given by $\mathcal{S}\triangleq \mathcal{S}^0\cup \mathcal{S}^*$. We consider binary investment variables for new thermal generators, i.e.,  $x^\mathcal{G}_g\in\{0,1\}$ for $g\in\mathcal{G}_T^*$, corresponding to each new installed unit.  With this formulation, we can represent the operating and commitment characteristics of the new thermal units. Since renewable generators and storage elements have more flexibility in how much capacity we can install, we consider continuous variables for these investments, i.e.,  $x^\mathcal{G}_g\in\mathbb{R}_+$ for $g\in\mathcal{G}_R^*$ and $x^\mathcal{S}_s\in\mathbb{R}_+$ for $s\in\mathcal{S}^*$.  All investment variables are scaled by a per unit capacity value (in MW) to yield the nameplate capacity of a new unit.  We also consider allowing for the retirement of a subset of the existing thermal generators, i.e.,  the set $\mathcal{G}_T^\dagger \subset \mathcal{G}_T^0$, and introduce the retirement decision variable $x^\dagger_g \in\{0,1\}$ for $g \in \mathcal{G}_T^\dagger$.  See Table \ref{tab:first_stage_variables} for a list of the first-stage variables and Table \ref{tab:sets_mappings} for a list of the sets and mapping functions used in the mathematical model .

We use the ``tight'' UC formulation introduced in \cite{knueven2020mixed}, a work which provides a comprehensive assessment of a variety of different UC model formulations.  This formulation uses a three-binary formulation to represent unit status for thermal generators, given in Table \ref{tab:second_stage_variables} by $u$, $v$, and $w$ variables.  It also introduces additional power variables for thermal generators.  See Table \ref{tab:second_stage_variables} for a complete list of second-stage variables.

\begin{table}[!t]
\caption{First-stage (investment) decision variables}
    \label{tab:first_stage_variables}
  \centering
      \rowcolors{2}{llnlLivermoriumIce}{white}
    \begin{tabular}{p{1.3cm} >{\raggedright\arraybackslash}p{3.4cm} p{2.5cm}}
    \rowcolor{llnlElementalNavy}
    \myth{Symbol} & \myth{Description} & \myth{Domain} \\    
    $x_g^{\mathcal{G}}$ & Per unit generator investment at candidate $g \in \mathcal{G}^*$ &
      \begin{tabular}[c]{@{}c@{}}$\{0,1\}$ if $g \in \mathcal{G}_T^*$\\ $\mathbb{R}_+$ if $g \in \mathcal{G}_R^*$\end{tabular} \\
    $x_s^{\mathcal{S}}$ & Per unit storage (power) investment at candidate $s \in \mathcal{S}^*$ & $\mathbb{R}_+$ \\
    $x_g^\dagger$ & Indicator for the retirement of existing generator $g \in \mathcal{G}_T^\dagger$ & $\{0,1\}$ \\
  \end{tabular}
\end{table}

\begin{table*}[!t]
\centering
\caption{Sets and mapping functions}
\label{tab:sets_mappings}
    \rowcolors{2}{llnlLivermoriumIce}{white}
\begin{tabular}{l l}
\rowcolor{llnlElementalNavy}
\myth{Symbol} & \myth{Description} \\
    $\mathcal{B}$ & Set of buses (nodes) \\
    $\mathcal{L}$ & Set of branches (transmission lines \& transformers) \\
    $\mathcal{L}^{\mathrm{f}}(b)$ & Set of branches whose ``from'' bus is $b$ \\
    $\mathcal{L}^{\mathrm{t}}(b)$ & Set of branches whose ``to'' bus is $b$ \\
    $\mathcal{G}$ & Set of all generators (existing and candidate) \\
    $\mathcal{G}^0$ & Set of existing generators \\
    $\mathcal{G}^0_T$ & Set of existing thermal generators \\
    $\mathcal{G}^0_R$ & Set of existing renewable generators \\
    $\mathcal{G}^*$ & Set of candidate generators \\
    $\mathcal{G}^*_T$ & Set of candidate thermal generators \\
    $\mathcal{G}^*_R$ & Set of candidate renewable generators \\
    $\mathcal{G}_T$ & Set of all thermal generators \\
    $\mathcal{G}_R$ & Set of all renewable generators \\
    $\mathcal{G}_T^\dagger$ & Set of existing thermal generators eligible for retirement \\
    $\mathcal{G}^0(b)$ & Set of existing generators located at bus $b$ \\
        $\mathcal{G}^*(b)$ & Set of candidate generators located at bus $b$ \\
    $\mathcal{G}^0(c,b)$ & Existing generators of technology type $c$ at bus $b$ \\
    $\mathcal{G}^*(c,b)$ & Candidate generators of technology type $c$ at bus $b$ \\
    $\mathcal{S}$ & Set of all storage elements (existing and candidate) \\
    $\mathcal{S}^0$ & Set of existing storage elements \\
    $\mathcal{S}^*$ & Set of candidate storage elements \\
    $\mathcal{S}^0(b)$ & Set of existing storage elements located at bus $b$ \\
        $\mathcal{S}^*(b)$ & Set of candidate storage elements located at bus $b$ \\
    $\mathcal{S}^0(c,b)$ & Existing storage elements of technology type $c$ at bus $b$ \\
    $\mathcal{S}^*(c,b)$ & Candidate storage elements of technology type $c$ at bus $b$ \\
    $\mathcal{C}^{\mathcal{G}}$ & Set of generation technology types (e.g., CCGT, PV, wind) \\
    $\mathcal{C}^{\mathcal{S}}$ & Set of storage technology types (e.g., pumped hydro) \\
    $\mathcal{F}$ & Set of fuel types \\
    $\mathcal{G}^{\mathrm{p}}(f)$ & Thermal generators using fuel $f$ as primary fuel \\
    $\mathcal{G}^{\mathrm{s}}(f)$ & Thermal generators using fuel $f$ as secondary fuel \\
    $\mathcal{H}_g$ & Set of start-up categories (e.g., hot/warm/cold) for thermal generator $g\in\mathcal{G}_T$ \\
    $\mathcal{T}$ & Set of time periods (e.g., hours within a representative day), given as $\{1,\dots,T\}$ \\
    $\Omega$ & Set of scenarios (representative days) \\
  \end{tabular}
\end{table*}

\begin{table}[!t]
\caption[Second-stage (operation) decision variables]{Second-stage (operation) decision variables.  All second-stage variables are defined for $(t,\omega) \in \mathcal{T} \times \Omega$ unless otherwise indicated.}
    \label{tab:second_stage_variables}
  \centering
      \rowcolors{2}{llnlLivermoriumIce}{white}
    \begin{tabular}{p{1.3cm} >{\raggedright\arraybackslash}p{5cm} p{1.3cm}}
    \rowcolor{llnlElementalNavy}
    \myth{Symbol} & \myth{Description} & \myth{Domain} \\    
    $p_{g,t,\omega}^{\mathcal{G}}$ & Power output of generator $g \in \mathcal{G}$ (MW) & $\mathbb{R}_+$ \\
    $\overline{p}_{g,t,\omega}^{\mathcal{G}}$ & Maximum available power (headroom) of $g \in \mathcal{G}_T$ (MW) & $\mathbb{R}_+$ \\
    $p_{s,t,\omega}^{\mathcal{S}}$ & Energy stored in storage element $s \in \mathcal{S}$ (MWh), defined for $\mathcal S\times\{0\}\cup\mathcal T\times \Omega$  & $\mathbb{R}_+$ \\
    $p_{s,t,\omega}^{\mathcal{S}+}$ & Charging power at storage element $s \in \mathcal{S}$ (MW) & $\mathbb{R}_+$ \\
    $p_{s,t,\omega}^{\mathcal{S}-}$ & Discharging power at storage element $s \in \mathcal{S}$ (MW) & $\mathbb{R}_+$ \\
    $p_{\ell,t,\omega}^{\mathcal{L}}$ & Power flow on transmission branch $\ell \in \mathcal{L}$ (MW) & $\mathbb{R}$ \\
    $p_{b,t,\omega}^{\text{shed}}$ & Load shed at bus $b\in\mathcal{B}$ (MW) & $\mathbb{R}_+$ \\
    $u_{g,t,\omega}$ & On/off status of thermal generator $g \in \mathcal{G}_T$, defined for $\mathcal{G}_T\times\{0\}\cup\mathcal T\times \Omega$  & $\{0,1\}$\\
    $v_{g,t,\omega}$ & Start-up indicator for thermal generator $g \in \mathcal{G}_T$ & $\{0,1\}$ \\
    $w_{g,t,\omega}$ & Shut-down indicator for thermal generator $g \in \mathcal{G}_T$ & $\{0,1\}$ \\
    $d_{s,t,\omega}$ & Discharging indicator for storage element $s\in\mathcal{S}$ & $\{0,1\}$ \\
    $e_{s,t,\omega}$ & Charging indicator for storage element $s\in\mathcal{S}$ & $\{0,1\}$ \\
    $u_{g,t,\omega}^{\mathrm{p}}$ & Indicator for $g \in \mathcal{G}_T$ on and using primary fuel & $\{0,1\}$ \\
    $u_{g,t,\omega}^{\mathrm{s}}$ & Indicator for $g \in \mathcal{G}_T$ on and using secondary fuel & $\{0,1\}$ \\
    $f_{g,t,\omega}$ & Total fuel consumed by generator $g \in \mathcal{G}_T$ (MMBtu) & $\mathbb{R}_+$ \\
    $f_{g,t,\omega}^{\mathrm{p}}$ & Primary fuel consumed by generator $g \in \mathcal{G}_T$ (MMBtu) & $\mathbb{R}_+$ \\
    $f_{g,t,\omega}^{\mathrm{s}}$ & Secondary fuel consumed by generator $g \in \mathcal{G}_T$ (MMBtu) & $\mathbb{R}_+$ \\
    $c_{g,t,\omega}^{\text{su}}$ & Start-up cost for thermal generator $g \in \mathcal{G}_T$ at time $t$ (\$/h) & $\mathbb{R}_+$ \\
    $ux_{g,t,\omega}$ & Lifted product variable approximating $u_{g,t,\omega} \, x_g^{\mathcal{G}}$ for $g \in \mathcal{G}_T$& $\{0,1\}$ \\
    $\sigma_{t,\omega}^{\text{r}}$ & Reserve shortfall (MW) & $\mathbb{R}_+$ \\
    $\sigma_{g,\omega}^{\text{util}}$ & Under-utilization of existing thermal generator $g \in \mathcal{G}_T^0$, across time periods in scenario (MWh), defined for $\mathcal{G}^0_T\times\Omega$ & $\mathbb{R}_+$
  \end{tabular}
\end{table}

\subsection{Objective}
The goal of CEP is to find the investment plan that minimizes investment-stage and operational-stage costs, given respectively as $C^\text{1}(\textbf{x})$ and $C^\text{2}(\textbf{x})$, where $\textbf{x}\in\mathbb{R}^{|\mathcal{G}^*|+|\mathcal{S}^*|+|\mathcal{G}_T^\dagger|}$ is some candidate investment and retirement plan.
Both cost terms are annualized to obtain the units of \$/year.
\begin{equation}
    \min_{\textbf{x}} C^\text{1}(\textbf{x}) +  T^{\text{rep}} \tau\cdot \sum_{\omega \in \Omega} \rho_\omega C_\omega^\text{2}(\textbf{x})\label{eqn:objective}
\end{equation}

The first-stage costs are given as the sum of annualized capital and annual fixed operations and maintenance (FOM) costs for new generators and storage elements as well as the costs to retire existing thermal generators.  We also include the FOM costs associated with existing generators and storage and assume that the retirement of a thermal generator eliminates its FOM costs.  The first-stage costs are a linear function of the first-stage investment and retirement decision variables:
\begin{multline}
C^\text{1}\triangleq
\sum_{g\in\mathcal{G}^*}
\overline{P}^{\mathcal G}_g
\left(C^{\mathcal G \text{-cap}}_{g} + C^{\mathcal G \text{-FOM}}_g\right)
x^{\mathcal G}_{g} \\
+ \sum_{s\in \mathcal{S}^*}
\overline{P}^{\mathcal S}_s
\left(C^{\mathcal S \text{-cap}}_{s}+C^{\mathcal S \text{-FOM}}_s\right)
x^{\mathcal S}_{s} \\
+ \sum_{g\in\mathcal{G}^\dagger_T}
C^{\mathcal{G}\text{-ret}}_g x^\dagger_g
+ C^{\mathcal G \text{-FOM}}_g \overline{P}^{\mathcal G}_{g}(1-x^\dagger_g) \\
+ \sum_{g\in \mathcal G^0\setminus \mathcal{G}^\dagger_T}
C^{\mathcal G \text{-FOM}}_g \overline{P}^{\mathcal G}_{g}
+ \sum_{s\in \mathcal{S}^0}
C^{\mathcal S \text{-FOM}}_s \overline{P}^{\mathcal S}_{s}
\label{eqn:invest_cost}
\end{multline}

The second-stage costs are given as the sum of generator production costs for thermal and renewable units, thermal generator commitment-stage costs,  storage costs, and penalties for load shedding, reserve shortfall, and thermal generator under-utilization:

\begin{align}
C_\omega^\text{2} \triangleq
& \sum_{t\in\mathcal T}\Bigg[\sum_{g\in \mathcal G_R} C^{\mathcal G\text{-VOM}}_g p^{\mathcal G}_{g,t,\omega}\notag\\
& \qquad + \sum_{g\in \mathcal G_T} \left( C^\text{prod}(p^{\mathcal G}_{g,t,\omega}) + c^\text{su}_{g,t,\omega} \right)\notag\\
& \qquad + \sum_{s\in \mathcal S} \left( C^{\mathcal S\text{-VOM}}_s + C^{\mathcal S-}_s \right) p^{\mathcal S-}_{s,t,\omega}
\Bigg]\notag
\\
& \quad + C^{\text{shed}} \sum_{t \in \mathcal T} \sum_{b\in \mathcal B} p^{\text{shed}}_{b,t,\omega}
+ C^\text{res}\sum_{t\in\mathcal{T}} \sigma^\text{r}_{t,\omega}
\notag\\
& \quad + C^\text{util}\sum_{g\in\mathcal{G}_T}\sigma^\text{util}_{g,\omega},
\quad \forall \omega\in\Omega\label{eqn:oper_cost}
\end{align}
where $C^\text{prod}(\cdot)$ corresponds to the thermal generator production costs as defined in Equation \eqref{eqn:prod_cost}. We assume that the operational horizon $\mathcal T$ considered has time intervals of uniform length $\tau$. There are $T^{\text{rep}}$ operational horizons per year.  The parameters introduced in \eqref{eqn:objective}, \eqref{eqn:invest_cost} and \eqref{eqn:oper_cost} are defined in Table \ref{tab:cost_parameters}.

\begin{table}[!t]
\caption{Parameters introduced in cost expressions \eqref{eqn:objective} through \eqref{eqn:oper_cost} and land-use constraints \eqref{eq:max_gen_bus} through \eqref{eq:max_bus2}}
\label{tab:cost_parameters}
  \centering
    \rowcolors{2}{llnlLivermoriumIce}{white}
  \begin{tabular}{p{1.3cm} p{6.3cm}}
  \rowcolor{llnlElementalNavy}
    \myth{Symbol} & \myth{Description} \\    
    $\tau$ & Duration of each time period (h) \\
    $T^{\mathrm{rep}}$ & Number of representative scenarios (e.g., days) per year \\
    $\rho_\omega$ & Probability or weight of scenario $\omega \in \Omega$ \\
    $\overline{P}_g^{\mathcal{G}}$ & Nameplate unit capacity (MW) for generator $g \in\mathcal{G}$ (MW per unit of $x_g^{\mathcal{G}}$ for $g\in\mathcal{G}^*$)\\    
    $\overline{P}_s^{\mathcal{S}}$ & Nameplate power inversion capacity (MW) for storage $s\in\mathcal{S}$ (MW per unit of $x_s^{\mathcal{S}}$ for $s\in\mathcal{S}^*$)\\    
    $C_g^{\mathcal{G}\text{-cap}}$ & Annualized capital cost per MW of candidate generator $g\in\mathcal{G}^*$ (\$/MWy)\\
    $C_g^{\mathcal{G}\text{-FOM}}$ & Annual FOM cost per MW of generator $g\in\mathcal{G}$ (\$/MWy)\\
    $C_s^{\mathcal{S}\text{-cap}}$ & Annualized capital cost per MW of candidate storage $s\in\mathcal{S}^*$ (\$/MWy)\\
    $C_s^{\mathcal{S}\text{-FOM}}$ & Annual FOM cost per MW of storage $s\in\mathcal{S}$ (\$/MWy)\\
    $C_g^{\mathcal{G}\text{-ret}}$ & Annualized cost to retire existing thermal generator $g\in\mathcal{G}^\dagger_T$ (\$/y)\\
    $C_g^{\mathcal{G}\text{-VOM}}$ & VOM cost for generator $g\in\mathcal{G}$ (\$/MWh)\\
    $C_s^{\mathcal{S}\text{-VOM}}$ & VOM cost for storage discharge at $s\in\mathcal{S}$ (\$/MWh) \\
    $C_s^{\mathcal{S}-}$ & Additional cost of storage discharge (e.g., degradation) at $s\in\mathcal{S}$  (\$/MWh)\\    
    $C^{\text{shed}}$ & Penalty cost of load shedding  (\$/MWh)\\
    $C^{\text{res}}$ & Penalty cost of reserve shortfall  (\$/MWh)\\
    $C^{\text{util}}$ & Penalty cost of thermal generator under-utilization  (\$/MWh)\\
    $P^{\mathcal{G}\text{-max}}_{c,b}$ & Maximum allowable generation (MW) of technology type $c$ at bus $b$\\
    $P^{\mathcal{G}\text{-max}}_b$ & Total maximum allowable generation (MW) at bus $b$\\
    $P^{\mathcal{S}\text{-max}}_{c,b}$ & Maximum allowable storage (MW) of technology type $c$ at bus $b$\\
    $P^{\mathcal{S}\text{-max}}_{b}$ & Total maximum allowable storage (MW) at bus $b$\\
  \end{tabular}
\end{table}


\subsection{Constraints on New Investments}

We limit the maximum potential generation capacity that can be installed of each technology type at each bus based on resource availability and land-use restrictions:

\begin{multline}
    \sum_{g\in \mathcal{G}^*(c,b)} \overline{P}^{\mathcal G}_gx^{\mathcal G}_{g} +\sum_{g\in \mathcal{G}^0(c,b)} \overline{P}^{\mathcal G}_{g} \le P^{\mathcal{G}\text{-max}}_{c,b},\\\forall ~ c \in \mathcal C^\mathcal{G}, b\in \mathcal B \label{eq:max_gen_bus}
\end{multline}

We also choose to limit the amount of new generation at each bus, irrespective of type, in order to better distribute the new generation around the island:
\begin{align}
  \sum_{g\in \mathcal{G}^*(b)} \overline{P}^{\mathcal G}_gx^{\mathcal G}_{g} +\sum_{g\in \mathcal{G}^0(b)} \overline{P}^{\mathcal G}_{g}  &\le P^{\mathcal{G}\text{-max}}_b,\enskip ~\forall ~ b \in\mathcal{B} \label{eq:max_bus1}
\end{align}

Similar constraints are defined for storage elements:
\begin{multline}
    \sum_{s\in \mathcal{S}^*(c,b)} \overline{P}^{\mathcal S}_sx^{\mathcal S}_{s} +\sum_{s\in \mathcal{S}^0(c,b)} \overline{P}^{\mathcal S}_{s} \le P^{\mathcal{S}\text{-max}}_{c,b},\\\forall ~ c \in \mathcal C^\mathcal{S}, b\in \mathcal B \label{eq:max_stor_bus}
\end{multline}
\begin{align}
  \sum_{s\in \mathcal{S}^*(b)} \overline{P}^{\mathcal S}_sx^{\mathcal S}_{s} +\sum_{s\in \mathcal{S}^0(b)} \overline{P}^{\mathcal S}_{s}  &\le P^{\mathcal{S}\text{-max}}_b,\quad ~\forall ~ b \in\mathcal{B} \label{eq:max_bus2}
\end{align}

The parameters introduced in \eqref{eq:max_gen_bus} through \eqref{eq:max_bus2} are defined in Table \ref{tab:cost_parameters}.

\subsection{Second-stage Constraints from UC Model}
In this section, we describe the constraints and costs for the high-fidelity representation of thermal generators as well as constraints for renewable generators, storage, transmission limits, power balance, reserves, and fuel modeling.  These constraints effectively describe a classical UC/ED problem.

\subsubsection{UC Constraints on Thermal Generators}
From the tight UC formulation in \cite{knueven2020mixed}, we have logical constraints that relate the three binary status variables:
\begin{multline}
     u_{g,t,\omega}-u_{g,t-1,\omega} = v_{g,t,\omega} - w_{g,t,\omega},\\\forall ~g\in\mathcal{G}_T, t\in\mathcal{T},\omega\in\Omega\label{eqn:Logical}
\end{multline}

We enforce that the state of a unit at the end of the time horizon is the same as its starting state at $t=0$:
\begin{equation}
   u_{g,0,\omega}= u_{g,T,\omega},\quad\forall ~g\in\mathcal{G}_T, \omega\in\Omega\label{eqn:UnitOnInitialRule}
\end{equation}

The power generated relates to the unit commitment variables by the constraints:
\begin{multline}
    \underline{P}_g^\mathcal{G}u_{g,t,\omega}\leq  p^\mathcal{G}_{g,t,\omega} \leq  \eta^\mathcal{G}_{g,t,\omega}\overline{P}_g^\mathcal{G}u_{g,t,\omega},\\\forall ~g\in\mathcal{G}_T,t\in\mathcal{T},\omega\in\Omega\label{eqn:EnforceGeneratorOutputLimitsPartB}
\end{multline}
where for units that are ``on,'' the power generation is lower bounded by the minimum power output and upper bounded by the nameplate capacity, scaled by the generator's availability factor.
    
The maximum power available at a given generator $g$ for reserve margin purposes, i.e.,  its headroom, is bounded by the nameplate capacity, scaled by the availability factor:
\begin{equation}
    \overline{p}^\mathcal{G}_{g,t,\omega} \leq  \eta^\mathcal{G}_{g,t,\omega}\overline{P}_g^\mathcal{G},\quad\forall g\in\mathcal{G}_T,t\in\mathcal{T},\omega\in\Omega
\end{equation}

We enforce up-time and down-time constraints on thermal generators via the constraints:
\begin{gather}
u_{g,t,\omega}-u_{g,t-1,\omega} \leq u_{g,t',\omega}, ~~\forall ~g\in\mathcal{G}_T, \omega\in\Omega, \nonumber\\
\quad\forall ~\{t,t'\in\mathcal{T}: t+1\leq t'\leq t+UT_g-1\}
\label{eqn:UpTime}\\
u_{g,t-1,\omega}-u_{g,t,\omega} \leq 1-u_{g,t',\omega}, ~~ \forall ~g\in\mathcal{G}_T, \omega\in\Omega, \nonumber\\
\quad\forall ~\{t,t'\in\mathcal{T}: t+1\leq t'\leq t+DT_g-1\},
   \label{eqn:DownTime}
\end{gather}


To enforce ramp up and down limits, we use the Damc\i-Kurt ramping constraints defined in \cite{knueven2020mixed}: 
\begin{align}
p^\mathcal{G}_{g,t,\omega}-p^\mathcal{G}_{g,t-1,\omega}
\leq\;&
(SU_g-\underline{P}_g^\mathcal{G}-RU_g)\,v_{g,t,\omega}
\nonumber\\
&+
(\underline{P}_g^\mathcal{G}+RU_g)\,u_{g,t,\omega}
-
\underline{P}_g^\mathcal{G}\,u_{g,t-1,\omega}\nonumber\\
&\qquad,~\forall g\in\mathcal{G}_T,t\in\{2,\dots,T\},\omega\in\Omega
\label{eqn:EnforceMaxAvailableRampUpRates}\\
p^\mathcal{G}_{g,t-1,\omega}-p^\mathcal{G}_{g,t,\omega}
\leq\;&
(SD_g-\underline{P}_g^\mathcal{G}-RD_g)\,w_{g,t,\omega}
\nonumber\\
&+
(\underline{P}_g^\mathcal{G}+RD_g)\,u_{g,t-1,\omega}
-
\underline{P}_g^\mathcal{G}\,u_{g,t,\omega}\nonumber\\
& \qquad,~\forall g\in\mathcal{G}_T,t\in\{2,\dots,T\},\omega\in\Omega
\label{eqn:EnforceScaledNominalRampDownLimits}
\end{align}

Note that we include several additional ramping constraints that provide additional tightening to better approximate the convex hull, as included in the tight formulation of \cite{knueven2020mixed}.  These details are omitted for clarity.
The parameters introduced in \eqref{eqn:Logical} through \eqref{eqn:EnforceScaledNominalRampDownLimits} are defined in Table~\ref{tab:ramp_su_sd_params}.

\begin{table}[!t]
\caption{Parameters introduced for modeling generator operation and costs}
\label{tab:ramp_su_sd_params}
  \centering
    \rowcolors{2}{llnlLivermoriumIce}{white}
\begin{tabular}{l p{6cm}}
\rowcolor{llnlElementalNavy}
    \myth{Symbol} & \myth{Description} \\    
    $T$& Final time in $\mathcal{T}$\\
    $\underline{P}_g^{\mathcal{G}}$  & Minimum power output (MW) of thermal generator $g\in\mathcal{G}_T$ when on \\
    $\eta_{g,t,\omega}^{\mathcal{G}}$ & Factor of nameplate capacity of generator $g\in\mathcal{G}$ available at period $t\in\mathcal{T}$ in scenario $\omega\in\Omega$ \\
    $UT_g$ & Minimum up time (hours) of thermal generator $g\in\mathcal{G}_T$ \\
    $DT_g$ & Minimum down time (hours) of thermal generator $g\in\mathcal{G}_T$ \\
    $RU_g$ & Hourly ramp up limit (MW/h) for thermal generator $g\in\mathcal{G}_T$ \\
    $RD_g$ & Hourly ramp down limit (MW/h) for thermal generator $g\in\mathcal{G}_T$ \\
    $SU_g$ & Start-up ramp limit (MW/h) for thermal generator $g\in\mathcal{G}_T$ \\
    $SD_g$ & Shut-down ramp limit (MW/h) for thermal generator $g\in\mathcal{G}_T$\\
     $\text{f}_g^\text{p}$ & Label for the primary fuel type of thermal generator $g\in\mathcal{G}_T$\\
    $\text{f}_g^\text{s}$ & Label for the secondary fuel type of thermal generator $g\in\mathcal{G}_T$\\
    $C^{\mathcal G \text{-fuel}}_{f}$  & Fuel cost per unit of fuel $f\in\mathcal{F}$ (\$/MMBtu)\\
    $\underline{T}_{h,g}$ & Minimum time offline (hours) corresponding to start-up category $h\in\mathcal{H}_g$ for thermal generator $g\in\mathcal{G}_T$\\
    $C^{SU}_{h,g}$ & Start-up cost (\$) corresponding to start-up category $h \in \mathcal{H}_g$ for thermal generator $g\in\mathcal{G}_T$\\
  \end{tabular}
\end{table}

\subsubsection{Production and Commitment Costs for Thermal Generators}
The production costs for thermal generators are given by a linear combination of variable operations and maintenance (VOM) and fuel costs:
\begin{multline}
    C^\text{prod}(p^{\mathcal{G}}_{g,t,\omega})=  C^{\mathcal G \text{-VOM}}_gp^{\mathcal{G}}_{g,t,\omega}+C^{\mathcal G \text{-fuel}}_{\text{f}_g^\text{p}} f^\text{p}_{g,t,\omega}\\+C^{\mathcal G \text{-fuel}}_{\text{f}_g^\text{s}} f^\text{s}_{g,t,\omega},\enskip\forall g\in\mathcal{G}_T, t\in\mathcal{T},\omega\in\Omega\label{eqn:prod_cost}
\end{multline}
where the fuel costs are split between primary and secondary fuel types.  The dual fuel modeling of thermal generators is described in more detail in Section \ref{sec:fuel_modeling}.

The startup costs for thermal generators are defined as:
\begin{align}
c^\text{su}_{g,t,\omega}
&\geq
C^{SU}_{h,g}
\left(
u_{g,t,\omega}
- \sum_{i=1}^{\underline{T}_{h,g}} u_{g,t-i,\omega}
\right),
\nonumber\\
&\qquad \forall~h\in\mathcal{H}_g, g\in\mathcal{G}_T, t\in\mathcal{T},\omega\in\Omega\label{eqn:startup_cost}
\end{align}
where this constraint applies for each start-up category $h\in\mathcal{H}_g$ (typically, hot/warm/cold).  Because the start-up costs $C^{SU}_{h,g}$ are taken to be increasing with increasing start-up times, e.g., start-up costs are higher for a warm startup than a hot one, this constraint will only be binding for the start-up category during which the generator is turned on.
Note that in our model, we use a mathematically equivalent version of this constraint, as described by the tight formulation in \cite{knueven2020mixed}, and refer the reader to that work for more details.  The parameters introduced in \eqref{eqn:prod_cost} and \eqref{eqn:startup_cost} are defined in Table \ref{tab:ramp_su_sd_params}.


\subsubsection{Renewable Generators and Storage}
For existing variable resources, we have maximum power constraints:
\begin{align}
p^{\mathcal G}_{g,t,\omega} \le \eta^{\mathcal G}_{g,t,\omega} \overline{P}^{\mathcal G}_g, ~\forall~ g\in \mathcal G_R^0, t\in \mathcal T, \omega \in \Omega \label{eqn:variable_gen_bound}
\end{align}

For existing storage elements, storage energy level is bounded by nameplate power inversion capacity multiplied by the storage element's energy-to-power ratio:
\begin{align}
p^{\mathcal S}_{s,t, \omega} & \le \lambda_s\overline{P}^\mathcal{S}_s, ~ \forall~ s \in \mathcal S^0, t \in \mathcal T, \omega \in \Omega\label{eqn:max_energy_storage}
\end{align}
where variations of these constraints for candidate projects are handled in Equations \eqref{eqn:ModifyMaximumPowerOutput} and \eqref{eqn:ModifyMaxStorageEnergy}.

For all storage elements, we have maximum power charge and discharge limits:
\begin{align}
p^{\mathcal S+}_{s,t, \omega} & \le \overline{P}^\mathcal{S}_se_{s,t,\omega}, ~ \forall~ s \in \mathcal S, t \in \mathcal T, \omega \in \Omega\label{eqn:EnforceStorageInputLimitsPartB}\\
    p^{\mathcal S-}_{s,t, \omega} & \le \overline{P}^\mathcal{S}_sd_{s,t,\omega}, ~ \forall~ s \in \mathcal S, t \in \mathcal T, \omega \in \Omega\label{eqn:EnforceStorageOutputLimitsPartB}
\end{align}

We also have constraints that describe the dynamics of storage charge/discharge:
\begin{multline}
p^{\mathcal S}_{s,t,\omega}
=
p^{\mathcal S}_{s,t-1,\omega}
+ \tau\left(
\eta^{\mathcal{S}+}_{s} p^{\mathcal S+}_{s,t,\omega}
- \frac{p^{\mathcal S-}_{s,t,\omega}}{\eta^{\mathcal S -}_s}
\right), \\
\forall ~s\in\mathcal{S}, t\in \mathcal{T},\omega\in\Omega
\label{eqn:EnergyConservation}
\end{multline}

We enforce a cyclic behavior with storage, i.e., we use the energy state of the last period of the time horizon as the starting state at $t=0$:
\begin{equation}
    p^{\mathcal S}_{s,0,\omega}=p^{\mathcal S}_{s,T,\omega},\quad\forall~s\in\mathcal{S}, \omega\in\Omega
\end{equation}

We enforce that charging and discharging cannot happen simultaneously via the constraints:
\begin{equation}
    e_{s,t,\omega} + d_{s,t,\omega} \leq 1,\quad \forall ~s\in\mathcal{S},t\in\mathcal{T},\omega\in\Omega \label{eqn:InputOutputComplementarity}
\end{equation}

The parameters introduced in \eqref{eqn:variable_gen_bound} through \eqref{eqn:InputOutputComplementarity} are defined in Table \ref{tab:renewable_storage_params}.

\begin{table}[!t]
\caption{Parameters introduced for modeling storage, transmission, and system operation.}
\label{tab:renewable_storage_params}
  \centering
    \rowcolors{2}{llnlLivermoriumIce}{white}
\begin{tabular}{p{1.5cm} p{6cm}}
\rowcolor{llnlElementalNavy}
    \myth{Symbol} & \myth{Description} \\    
    $\lambda_s$ & Energy-to-power ratio of storage element $s\in\mathcal{S}$ (hours)\\
    $\eta_s^{\mathcal{S}+}$ & Charging efficiency of storage element $s\in\mathcal{S}$\\
    $\eta_s^{\mathcal{S}-}$ & Discharging efficiency of storage element $s\in\mathcal{S}$\\
    $D_{b,t,\omega}$ & Demand (MW) at bus $b$ at time $t$ in scenario $\omega$\\
    $\overline{P}_\ell^{\mathcal{L}}$ & Thermal limit (MW) of transmission branch $\ell\in\mathcal{L}$\\
    $RM_t$& Reserve margin requirement (MW) at time $t$\\
  \end{tabular}
\end{table}

\subsubsection{Transmission \& System Constraints}
To enforce nodal power balance from Kirchhoff's laws, we use the constraint:
\begin{multline}
\sum_{g \in \mathcal G(b)} p^{\mathcal G}_{g,t,\omega}
+ \sum_{s \in \mathcal S(b)}
\left(p^{\mathcal S-}_{s,t,\omega} - p^{\mathcal S+}_{s,t,\omega}\right) \\
+ \sum_{\ell \in \mathcal{L}^\text{t}(b)} p^{\mathcal L}_{\ell,t,\omega}
- \sum_{\ell \in \mathcal{L}^\text{f}(b)} p^{\mathcal L}_{\ell,t,\omega}
+ p^{\text{shed}}_{b,t,\omega}
= D_{b,t,\omega}, \\
\forall ~ b\in \mathcal B, t\in \mathcal T, \omega \in \Omega
\label{eq:PowerBalance}
\end{multline}

Constraints on branch flow thermal limits are given as:
\begin{equation}
    |p^\mathcal{L}_{\ell,t,\omega}| \leq \overline{P}^\mathcal{L}_\ell,\quad \forall \ell\in\mathcal{L},t\in \mathcal{T},\omega\in\Omega\label{eqn:LinePowerUB}
\end{equation}

We limit maximum load shed to demand:
\begin{align}
p^{\text{shed}}_{b,t,\omega} \le D_{b,t,\omega},\quad \forall~ b \in \mathcal B, t \in \mathcal T, \omega \in \Omega
\end{align}

The reserve shortfall is limited by the reserve margin requirement:
\begin{align}
    \sigma^\text{r}_{t,\omega} \leq RM_t,\quad\forall t\in\mathcal{T},~\omega\in\Omega
\end{align}

And the reserve shortfall $\sigma^{\text{r}}_{t,\omega}$ is determined by the constraint:
\begin{multline}
\sum_{b\in\mathcal{B}} D_{b,t,\omega} + RM_t \leq
\sum_{g\in\mathcal{G}_T} \overline{p}^\mathcal{G}_{g,t,\omega}
+ \sum_{g\in\mathcal{G}_R} p^\mathcal{G}_{g,t,\omega} \\
+ \sum_{s\in\mathcal{S}}
\left(p^{\mathcal{S}-}_{s,t,\omega} - p^{\mathcal{S}+}_{s,t,\omega}\right)
+ \sum_{b\in\mathcal{B}} p^\text{shed}_{b,t,\omega}
+ \sigma^\text{r}_{t,\omega}, \\
\forall t\in\mathcal{T},\omega\in\Omega
\label{eqn:EnforceReserveRequirements}
\end{multline}

The parameters introduced in \eqref{eq:PowerBalance} through \eqref{eqn:EnforceReserveRequirements} are defined in Table \ref{tab:renewable_storage_params}.


\subsubsection{Fuel Modeling}\label{sec:fuel_modeling}
The fuel consumed by a thermal generator relates to its energy production via the expression:
\begin{equation}
    f_{g,t,\omega} = HR_g(\tau p^{\mathcal{G}}_{g,t,\omega}),\quad\forall ~ g\in\mathcal{G}_T,t\in \mathcal{T},\omega\in\Omega\label{eqn:heat_rate_equation}
\end{equation}
where the heat rate function $HR_g(\cdot)$ may be either a linear or a convex linear-piecewise function of energy generated.  See \cite{knueven2020mixed} for the complete model of linear-piecewise cost curves.

For some thermal generators, we consider multiple fuel types, with the ability to switch between primary and secondary types.  For concision, we write out the equations below as if all thermal generators have dual fuel capabilities, but the equations are generalizable to the case where only a subset of thermal generators have this capability.  Note that this formulation is also generalizable to the case where we have location-specific fuels, as we do in the Puerto Rico test system, by indexing the different fuel types by location in addition to type.

The sum of all fuel used at a generator must be the sum of primary and secondary fuels used at that generator:
\begin{equation}
    f_{g,t,\omega} = f^\text{p}_{g,t,\omega} + f^\text{s}_{g,t,\omega},\quad \forall g\in\mathcal{G}_T,t\in\mathcal{T},\omega\in\Omega
\end{equation}
 
The constraints that relate the unit status variables with dual fuel status variables are given as:
\begin{equation}
    u_{g,t,\omega} = u^\text{p}_{g,t,\omega} +u^\text{s}_{g,t,\omega},\quad\forall g\in\mathcal{G}_T,t\in\mathcal{T},\omega\in\Omega\label{eqn:UnitOnLink}
\end{equation}

Where these status variables relate to fuel consumed variables via the constraints:
\begin{align}
    f^\text{p}_{g,t,\omega} \leq  \overline{F}_{g}u^\text{p}_{g,t,\omega},\quad\forall g\in\mathcal{G}_T,t\in\mathcal{T},\omega\in\Omega \label{eqn:PrimaryFuelConstr}\\
    f^\text{s}_{g,t,\omega} \leq  \overline{F}_{g}u^\text{s}_{g,t,\omega},\quad\forall g\in\mathcal{G}_T,t\in\mathcal{T},\omega\in\Omega \label{eqn:AuxiliaryFuelConstr}
\end{align}

Note that our formulation also includes similar fuel switching variables for unit start and stop and a related logical constraint like that in \eqref{eqn:Logical}.  These are omitted here for concision.

A limited fuel supply is considered for each fuel $f$ and scenario $\omega$ and is enforced for the total fuel consumed over the operational scenario:
\begin{multline}
\sum_{t\in\mathcal{T}}
\left(
\sum_{g\in\mathcal{G}^\text{p}(f)} f^\text{p}_{g,t,\omega}
+ \sum_{g\in\mathcal{G}^\text{s}(f)} f^\text{s}_{g,t,\omega}
\right)
\leq FS_{f,\omega}, \\
\forall~f\in\mathcal{F},\omega\in\Omega
\label{eqn:FuelLimitConstr}
\end{multline}

The parameters introduced in \eqref{eqn:heat_rate_equation} through \eqref{eqn:FuelLimitConstr} are defined in Table \ref{tab:fuel_params}.

\begin{table}[!t]
\caption{Parameters introduced for fuel modeling}
\label{tab:fuel_params}
  \centering
    \rowcolors{2}{llnlLivermoriumIce}{white}
\begin{tabular}{p{1.5cm} p{6cm}}
\rowcolor{llnlElementalNavy}
    \myth{Symbol} & \myth{Description} \\    
    $HR_g(\cdot)$ & Heat rate function (fuel input as a function of energy output) for thermal generator $g\in\mathcal{G}_T$ (MMBtu/MWh)\\
    $\overline{F}_g$ & Maximum fuel use rate for thermal generator $g\in\mathcal{G}_T$ (MMBtu/h) \\
    $FS_{f,\omega}$ & Total fuel supply over the scenario horizon of fuel $f\in\mathcal{F}$ in scenario $\omega\in\Omega$ (MMBTu/day)\\
  \end{tabular}
\end{table}

\subsection{Constraints Linking First- \\ and Second-stage Decisions }

The amount of power generated at a candidate generator cannot exceed its installed capacity, adjusted by the fraction of capacity that is available for the given scenario and time period.  For new generators, we have:
\begin{align}
p^{\mathcal G}_{g,t,\omega} \le \eta^{\mathcal G}_{g,t,\omega} \overline{P}^{\mathcal G}_gx^{\mathcal G}_{g}, \quad\forall~ g\in \mathcal G^*, t\in \mathcal T, \omega \in \Omega \label{eqn:ModifyMaximumPowerOutput}
\end{align}
where for new variable resources, such as wind and solar power, $\eta^{\mathcal G}_{g,t,\omega}$ will vary in time depending on the availability of the resource. For new thermal generators, $\eta^{\mathcal G}_{g,t,\omega}$ is the derate factor through which we model forced and planned outages.

For thermal generators, we also enforce the following constraint to ensure the headroom variable reflects whether or not a candidate generator is built:
\begin{equation}
    \overline{p}^{\mathcal G}_{g,t,\omega} \leq \eta^{\mathcal G}_{g,t,\omega} \overline{P}^{\mathcal G}_gx^{\mathcal G}_{g},  \quad\forall~ g\in \mathcal G_T^*, t\in \mathcal T, \omega \in \Omega
\end{equation}

For existing thermal generators that retire, we enforce that units must be off in the case that they retire:
\begin{align}
u_{g,t,\omega}\leq 1-x^\dagger_g,\quad~\forall~g\in\mathcal{G}_T^\dagger,t\in\mathcal{T},\omega\in\Omega
\end{align}

For new thermal generators, we also add lower bounds on minimum power output:
\begin{align}
    p^{\mathcal G}_{g,t,\omega} \geq  \underline{P}^{\mathcal G}_{g} ux_{g,t,\omega}, \quad~\forall~ g\in \mathcal G_T^*, t\in \mathcal T, \omega \in \Omega \label{eqn:ModifyMinimumPowerOutput}
\end{align}
where the variable $ux_{g,t,\omega}$ is the lifted binary variable corresponding to $u_{g,t,\omega}\cdot x^\mathcal{G}_{g}$ and constrained by the linear McCormick envelope:
\begin{align}
    ux_{g,t,\omega} &\leq \min\{x^\mathcal{G}_{g},u_{g,t,\omega}\},\quad \forall g\in\mathcal{G}_T^*,t\in \mathcal T, \omega \in \Omega\\
ux_{g,t,\omega} & \geq x^\mathcal{G}_{g} +  u_{g,t,\omega} - 1,\quad \forall g\in\mathcal{G}_T^*,t\in \mathcal T, \omega \in \Omega
\end{align} 

Similarly, for new storage elements, we have limits on energy storage and power charging and discharging, based on the installed capacity:
\begin{align}
p^{\mathcal S}_{s,t, \omega} & \le \lambda_s\overline{P}^\mathcal{S}_sx^{\mathcal S}_{s}, \quad\forall~ s \in \mathcal S^*, t \in \mathcal T, \omega \in \Omega\label{eqn:ModifyMaxStorageEnergy}\\
p^{\mathcal S-}_{s,t, \omega} & \le \overline{P}^\mathcal{S}_sx^{\mathcal S}_{s}, \quad\forall~ s \in \mathcal S^*, t \in \mathcal T, \omega \in \Omega\\
p^{\mathcal S+}_{s,t, \omega} & \le \overline{P}^\mathcal{S}_sx^{\mathcal S}_{s}, \quad \forall~ s \in \mathcal S^*, t \in \mathcal T, \omega \in \Omega
\end{align}

To discourage the retention of costly thermal generators that are rarely dispatched, we introduce a soft minimum-utilization requirement for existing thermal units:
\begin{align}
\sum_{t\in\mathcal{T}} p^{\mathcal{G}}_{g,t,\omega}
+ \sigma^\text{util}_{g,\omega}
&\geq
\kappa_g \overline{P}^\mathcal{G}_g |\mathcal{T}| (1-x^\dagger_g), \enskip\forall g\in\mathcal{G}^\dagger_T, \omega\in\Omega\\
\sum_{t\in\mathcal{T}} p^{\mathcal{G}}_{g,t,\omega}
+ \sigma^\text{util}_{g,\omega}
&\geq
\kappa_g \overline{P}^\mathcal{G}_g |\mathcal{T}|, \enskip\forall g\in\mathcal{G}^0_T\setminus \mathcal{G}^\dagger_T, \omega\in\Omega
\label{eq:min_gen_utilization}
\end{align}
where $\kappa_g$ is the minimum utilization factor for existing thermal generator $g$. Utilization below that threshold is penalized linearly via the violation term $\sigma^\text{util}_{g,\omega}$ at cost $C^{\text{util}}$.   Note that this constraint is intended as a proxy for economic pressure to retire underused units, not as a physical operating requirement.

\subsection{Assumptions}
The formulation described above includes the following assumptions:
\begin{itemize}
\item The sets of buses and transmission branches are assumed to be fixed, i.e., we do not consider bus additions or transmission expansion in this work.  Transmission expansion is unlikely during the time horizon of the study.
\item The distribution system is not considered in the modeling, so costs associated with distribution system improvements are not included.
\item The model is technology-agnostic, i.e., policy concerns like renewable portfolio standards (RPS) are not included. Sensitivities with renewables-only future expansion are included and explicitly noted; however, only some of these correspond to having no existing thermal units in the future system (a surrogate for a 100\% RPS requirement).
\item All new generators are assumed to be utility-scale and connected to the transmission system.  We do not consider DERs as distinct from behind-the-meter load.
\item Storage and renewable generator investments are treated as continuous decisions in our model.
\item We only consider instantaneous generation and storage expansion decisions, i.e.,  investments selected by the model are available immediately for operation.
\item The model considers all operational hours of the high-stress 3-month time horizon from May to July, unless otherwise noted.
\item All costs are scaled to be in units of \$/year.  Capital costs are amortized over the lifetime of the plant during data processing, and operational costs are calculated over a representative 3-month time horizon and scaled to a year.
\item We neglect real-time and day-ahead forecast errors in the commitment and operation of the generation dispatch; however, we choose a large reserve margin that helps compensate for this.
\item Generator capacities and costs are taken to be after generation losses.
\item Operations and maintenance costs for transmission lines and transformers are assumed to be constant across all cases so are removed for simplicity.
\item The transmission system is modeled with a nodal ``pipe-and-bubble'' approach, which treats power flow like a transportation network with branches subject to thermal limits.  This approach neglects the physics of power flow and overestimates the transfer capacity in the system.  This assumption was made for computational reasons, rather than modeling needs, and could be improved in a future iteration of this work.
\item We do not explicitly model transmission outages or an $N-1$ criterion. Outages are addressed via the reserve margin constraint, which we take to be approximately twice the size of the largest thermal unit in the existing system.
\item UC decisions are made over representative days.  Consecutive commitment and operational decisions are represented within a day but not across different days.  The ending on/off status of a thermal generator at the last time period is taken to be the initial status of the generator at the start of the same day.  However, minimum up/down time and ramping are not enforced for this end-of-day/start-of-day linking.
\item Similarly for storage modeling, we assume that the ending storage level at the last time period is the same as the initial storage level at the start of the same day.  This models charge/discharge cycles.  However, a limitation of this approach is that we do not capture the value of storage across days.
\item We assume that retiring a generator eliminates FOM costs for that generator.  Thermal units eligible for retirement will be retired only if doing so saves on overall system costs, as determined by the retirement costs, FOM cost savings, generator under-utilization penalty, and system operation.  
\item The under-utilization of existing thermal generators is computed for each day (summing over 24 hours) in the 3-month horizon, with soft enforcement via a penalty term.  This formulation allows for model decomposition via the progressive hedging algorithm.
\item Forced and planned outages are modeled for new thermal generators via derating to meet an assumed availability factor, rather than as a binary time-series input as we do for existing thermal generators.
\end{itemize}

\subsection{Limitations of the Approach}
There are several limitations of our approach, namely that we are not capturing power flow (mentioned in the section above) but also that we are not capturing sub-hourly variability and system dynamics/stability.  Because of this, there could be improved investments for an objective based on system stability, e.g., Battery Energy Storage Systems (BESS) may have more value in that context.  Post-processing can be used to validate power flow and reserve margins on a set of pre-specified contingencies, i.e.,  branch or generator outages, and our colleagues at PNNL have been validating our results using the PSS/E software package. We note that nearly all CEP models share these and other limitations, due to computational considerations.
Because our model is only a two-stage model, we do not capture multi-stage investment approaches that may be preferred for system reliability, i.e.,  gradual or sequential buildout approaches.  Certain resources such as BESS may be more costly long-term but quicker to deploy due to supply chain constraints, regulatory approvals, and construction time.  These facts are not within the scope of our model but we do consider a sensitivity on the types of generation available for near-term buildouts in Section \ref{sec:gradual_buildout}.  Finally, the presented study is a ``snapshot'' run for the season of highest stress in the system and should not be viewed as a full-scale resource adequacy study.  For a complete resource adequacy study, we would need to stress test the system with a variety of different weather, demand, forced outage, and fuel price scenarios, while considering an N-1 contingency criterion.  This was not within scope for this work.

In summary, the presented results are best interpreted as identifying promising generation and storage portfolios for further evaluation, rather than as final implementation recommendations. Candidate portfolios should be supplemented by resource adequacy analysis, production cost validation across additional weather years, transmission security assessment, dynamic stability studies, and implementation feasibility review.

\section{Puerto Rico (PR) Test System}
\label{sec:test_system}
\subsection{Existing System}

To obtain a list of existing and planned generators and storage elements, various publicly available sources were consulted and compared. Most basic characteristics like nameplate capacities, generation technology, compatible fuel types and heat rates by unit or power plant are available in the monthly and daily generation reports by LUMA \cite{luma_monthly_2025,luma_daily_gen_report} and their most recent resource adequacy studies \cite{luma_ra_study_2025}.
A series of detailed technical reports performed by Sargent \& Lundy are publicly available for most of the Genera thermal units and were also consulted for validating and verifying nameplate, configuration, and status data \cite{sl_ie_report_1,sl_ie_report_2,sl_ie_report_3,sl_ie_report_4,sl_ie_report_5,sl_ie_report_6,sl_ie_report_7}.
Detailed operational parameters for the Unit Commitment model for all thermal units were obtained from a PLEXOS Production Cost Model (PCM) provided by PNNL.  An overview of the data sources used for each technoeconomic parameter is provided in Table \ref{tab:technoeconomic_existing_system}.

\begin{table*}[ht]
\centering
\caption{Data sources for technoeconomic parameters for existing system}
\label{tab:technoeconomic_existing_system}
    \rowcolors{2}{llnlLivermoriumIce}{white}
\begin{tabular}{p{6cm} p{10cm}}
\rowcolor{llnlElementalNavy}
\myth{Data item} & \myth{Source} \\
Breakout of power plants into individual units &
LUMA daily generation reports \cite{luma_daily_gen_report}, validated by Sargent \& Lundy reports \cite{sl_ie_report_1,sl_ie_report_2,sl_ie_report_3,sl_ie_report_4,sl_ie_report_5,sl_ie_report_6,sl_ie_report_7} \\
Status of units &
LUMA daily generation reports; 2025 planned outages from \cite{luma_ra_study_2025} \\
Nameplate capacity (MW) for units &
``Dependable capacity'' in LUMA Resource Adequacy Study \cite{luma_ra_study_2025} \\
Minimum stable level (MW) for units &
Sargent \& Lundy reports \cite{sl_ie_report_1,sl_ie_report_2,sl_ie_report_3,sl_ie_report_4,sl_ie_report_5,sl_ie_report_6,sl_ie_report_7} \\
Fuel compatibility for generation units &
LUMA RA Study \cite{luma_ra_study_2025}, validated by Sargent \& Lundy reports \cite{sl_ie_report_1,sl_ie_report_2,sl_ie_report_3,sl_ie_report_4,sl_ie_report_5,sl_ie_report_6,sl_ie_report_7} and conversations with DOE sponsor \\
Plant heat rates &
December 2024 monthly generation report by LUMA \cite{luma_monthly_2025} \\
Ramp up/down rates &
PLEXOS PCM provided by PNNL \\
Minimum up/down times &
PLEXOS PCM provided by PNNL \\
FOM costs for existing units &
Sargent \& Lundy reports \cite{sl_ie_report_1,sl_ie_report_2,sl_ie_report_3,sl_ie_report_4,sl_ie_report_5,sl_ie_report_6,sl_ie_report_7} \\
VOM costs for existing units &
Sargent \& Lundy reports \cite{sl_ie_report_1,sl_ie_report_2,sl_ie_report_3,sl_ie_report_4,sl_ie_report_5,sl_ie_report_6,sl_ie_report_7} \\
Cost conversion to 2024 USD &
Bureau of Labor Statistics (BLS) WPU10 Producer Price Index for metals and metal products \cite{bureau_of_labor_statistics_wpu10_2025} \\
Fuel supplies of coal, NG, bunker, diesel &
Annual Electricity Fuel Consumption (MMBtu) from May--July 2024 EIA-923 form \cite{us_energy_information_administration_form_2023} \\
\end{tabular}
\end{table*}

Power flow models of Puerto Rico's transmission system at voltage levels of 38kV and above were obtained from the team at PNNL in PSS/E RAW format for the current system (2025) and the post-stabilization planned (2027-2030) system.
A manual mapping of units in the power flow cases to existing and planned units was performed to identify the connection point of each generation unit. When mismatches between the operational status or nameplate capacity of units in the power flow models and the rest of the data gathered by the LLNL team were found, the power flow case data were superseded.
The load data in the power flow cases were used as a baseline to perform a system-wide normalization and obtain load distribution factors to apply to the load time-series described in Section \ref{sec:time-series}.
Radial portions of the system with voltage levels below 38kV were aggregated to the closest substation at or above 38kV.  The status of transmission branches and loads was not modified from the values found in the power flow cases. Transmission branches and load elements set as ``out-of-service'' are assumed to be non-operational due, for example, to still unrepaired prior hurricane damage.

\subsection{Data for New Generators and Storage Elements}

The list of planned generators and storage elements was obtained from the most recent stabilization plan \cite{nepr_stab_plan_2025}.
The statuses of many of these projects have been dynamic during recent months and were updated during execution of the study after conversations with DOE sponsors and the Puerto Rico Public-Private Partnerships Authority (P3A), most recently in Q3 of 2025.  These planned generators are taken as an exogenous input to the model and provide a baseline system on top of which new generation is selected via CEP. 
We consider various sensitivities to the assumptions around which planned generators will materialize in Section \ref{sec:cep_results}.

For new generation considered as part of the CEP model, we looked at 9 different generation types and one type of storage (i.e., BESS).  The list of generation types considered can be found in Table \ref{tab:gentypes}.
Costs and technical parameters for thermal units were provided by the DOE Office of Fossil Energy (DOE-FE), now the Office of Hydrocarbons and Geothermal Energy.  DOE-FE based these on the National Energy Technology Laboratory (NETL) report \cite{netl_report}, where these costs were scaled by capex-specific escalators based on Chemical Engineering Plant Cost Index (CEPCI) data and location-specific PR escalators from DOE. Costs and parameters for other technologies were obtained from the Annual Energy Outlook published by the U.S. Energy Information Administration (EIA), applying some locational overhead to Puerto Rico. Values are provided in Table \ref{tab:costs_new_gen}.

\begin{table*}[ht]
\centering
\caption{List of generation types considered in the study}
\label{tab:gentypes}
\begin{tabular}{p{7cm} p{7cm}}
\rowcolor{llnlElementalNavy}
\multicolumn{2}{c}{\myth{Existing and not considered for new units}} \\
\rowcolor{white} Residual-Fuel-Oil-powered Steam Turbine & Coal-fired Steam Turbine \\
\rowcolor{gray!15}
Hydropower & Landfill gas \\
\rowcolor{llnlImpactBlue}
\multicolumn{2}{c}{\myth{In existing system, also considered for new units}} \\
\rowcolor{white} Utility-scale Solar PV & Onshore Wind \\
\rowcolor{llnlEnergeticAzure}
\multicolumn{2}{c}{\myth{Not in existing system but considered for new units}}\\

\rowcolor{white} 1x1 Combined Cycle (H-Class) & 1x1 Combined Cycle (F-Class) \\
\rowcolor{gray!15}
Large Simple Cycle (F-Class) & Small Simple Cycle (1 large aero) \\
\rowcolor{white}Small Simple Cycle (2 small aero) & Small Simple Cycle (6 large RICE) \\
\rowcolor{gray!15}
Small Simple Cycle (12 small RICE) & \\
\end{tabular}
\end{table*}

\subsection{Time-series Data for Load and \\Generator Availability}
\label{sec:time-series}
The load time-series data for all scenarios are the historical, hourly load time-series data provided by LUMA. Historical data from LUMA were provided from 2022-2024, with 2024 data used in all scenarios described in this report (future loads are uniformly scaled variations of the 2024 load). Load time-series data in all scenarios correspond to the load served by LUMA and thus do not account for any unserved load, which may underestimate the true demand on the system.  The aggregate load is distributed to in-service buses in the system based on the load values in the PSS/E test system. For baseline simulations on the existing system, historical served load from 2024 is used. For future expansion scenarios corresponding to the year 2030, the historical 2024 load is increased by 7\% overall, across all hours and buses, where this 7\% is derived from estimates in previous Resource Adequacy studies \cite{luma_ra_study_2025}.  We consider two additional sensitivities to future load, one with elevated loads and one with reduced loads, compared to this 2030 ``base case.''
 
Wind and solar generation time-series data used in this analysis are also historical, hourly records provided by LUMA. The data provided were the power generation for each utility-scale wind and solar generation plant in the system. We use the nameplate capacity for each plant in the PLEXOS PCM to convert these generation time-series data (in MW) to capacity factors (0-1), which are used as the generation availability. We compared these calculated capacity factors (based on LUMA actual data) to capacity factors from the PLEXOS PCM (see Figure \ref{fig:gen_avail} in the Appendix). In some cases, the offset between both datasets is relatively small (3.7 to 7.2\% at Oriana and Fonoroche Humacao, respectively), while in other cases the PLEXOS PCM data significantly overestimate 2024 capacity factors compared to historical records from LUMA (61.6\% at San Fermin). Given this difference, we use calculated capacity factors from LUMA actuals for all scenarios. For future expansion scenarios, we use the same 2024 weather year, but note that future analyses could include additional weather years not considered here.

\subsection{Forced Outages}

The high rate of forced outages for existing generation is the main culprit for the low quality of service currently experienced by Puerto Rico customers.
To accurately capture this reality in our study, forced outages were incorporated in the operational scenarios generated for the analysis.
Forced outages for each existing plant were randomly generated such that the fraction of hours out of service matches the outage rate specified in \cite{luma_ra_study_2025}, and each outage was assumed to last 40 hours, consistent with \cite{luma_ra_study_2025}.
Planned outages included in Table 12 of \cite{luma_ra_study_2025} were also incorporated into the generated outage scenarios.  For new generation units, no individual outages are generated. We instead represent any required maintenance or forced outage via an 85\% capacity factor, which we model by derating the nameplate capacity of all new units by this amount at all periods.
No transmission outages were explicitly modeled.

\subsection{Fuel Supplies and Costs}\label{sec:fuel_supply_data}

The fuels considered in the PR test system are natural gas (NG), bunker (also called residual fuel oil -- RFO), diesel, and coal.  We used the EIA-923 form to obtain data on Puerto Rico's fuel consumption for electricity production from May to July 2024 \cite{us_energy_information_administration_form_2023}.  This plant-specific data was used to obtain the fuel consumption at the locations of all  existing thermal plants.  Because this is historical \textit{consumption} data, it likely underestimates the possible available \textit{supply} of fuels.  For example, historical consumption of fuel may have been reduced if some generators were out-of-service, but that does not necessarily imply a lack of available imports.  After several rounds of fuel supply sensitivity studies and discussions with DOE staff, a set of supply scenarios was convened; it is described in Table~\ref{tab:fuel_scenarios}.
%
\begin{table*}[ht]
\centering
\caption{Fuel supply scenarios}
\label{tab:fuel_scenarios}
\begin{tabular}{p{4cm} p{10cm}}
\rowcolor{llnlElementalNavy}
\myth{Region}&\myth{Existing (2024) fuel supply} \\
\rowcolor{white}
San Juan (SJ)&fix NG and bunker to the 2024 consumption levels from EIA-923\\
    \hline
    \rowcolor{gray!15}
    Palo Seco (PS)&allow 20\% of the 2024 NG consumption at SJ (can be trucked to PS), fix bunker to the 2024 consumption level\\
    \hline
    \rowcolor{white}
    Aguirre& fix bunker to the 2024 consumption level\\
    \hline
    \rowcolor{gray!15}
     Costa Sur (CS)& allow unlimited NG and bunker\\
     \hline
     \rowcolor{white}
     EcoElectrica& allow unlimited NG\\
     \hline
     \rowcolor{gray!15}
     AES plant& fix coal to the 2024 consumption level\\
     \hline
     \rowcolor{white}
   All locations except AES & Allow unlimited diesel -- since it is the most expensive fuel, it will only be used when other fuels are not available \\
\rowcolor{llnlImpactBlue}
\myth{Region}&\myth{Future (2030) fuel supply} \\
\rowcolor{white} 
San Juan, Palo Seco, Cambalache, and Mayaguez & Unlimited NG \\
\rowcolor{gray!15}
All others & Same as existing \\
\rowcolor{llnlEnergeticAzure}
\myth{Region}&\myth{Future+ expanded fuel supply}\\
\rowcolor{white} 
Aguirre \& AES & Unlimited NG \\
\rowcolor{gray!15}
All others & Same as futures case \\
\end{tabular}
\end{table*}

Note that the fuel supply for each representative day in Equation (\ref{eqn:FuelLimitConstr}) is taken as the daily average consumption (in MMBtu/day) from the May to July 2024 EIA data.  This may be overly conservative as it does not allow for shifting fuel from high to low load days.

The location-specific fuel costs for San Juan, Palo Seco, Costa Sur, EcoElectrica, and the new Energiza project were based on DOE input for NG prices, which are relative to Henry Hub (HH) prices.  We use historical HH prices from May to July 2024 for the existing cases and projected HH prices for 2030 (based on the reference case in EIA's 2025 AEO) for the future cases \cite{EIA_HenryHub_Monthly,us_energy_information_administration_annual_2025}.  We scale the 2030 price by a seasonality factor for the May to July period, based on the 2024 historical HH data.  The location-specific bunker fuel prices for San Juan, Palo Seco, Aguirre, and Costa Sur are obtained from the PREPA PI dashboard, from a snapshot of prices taken in October 2025.  The diesel fuel prices were provided by DOE based on present-day oil prices and are not location-specific. For the remaining fuel costs, we used the LUMA December 2024 report, which has prices for NG, bunker, and coal but does not split these out by location within Puerto Rico \cite{luma_monthly_2025}. See Section \ref{sec:fuel_costs} in the Appendix for more details on the actual fuel costs used in the model.

\subsection{Additional Data Assumptions}
We also make the following assumptions to fill in unknown or missing data.
\begin{itemize}
\item We assume land constraints at existing wind and solar sites allow for expansion up to 2x existing capacity, unless otherwise noted. 
\item The Value of Lost Load (VOLL), $C^{\text{shed}}$ in our model, is taken to be \$30,000/MWh.
\item The reserve shortfall penalty $C^{\text{res}}$ is taken to be \$2,000/MWh.
\item The generator under-utilization penalty $C^{\text{util}}$ is taken to be \$2,000/MWh. 
\item We do not explicitly model growth in behind-the-meter generation.
\item The reserve margin target is taken to be 900 MW, about 2x the size of the largest thermal unit in the existing system.
\item We only consider retiring existing thermal generators in the system, not existing renewable generators or new builds.  We do not allow retirements of units with existing contractual agreements, such as AES and EcoElectrica.
\item We use the default minimum capacity factors ($\kappa_g$ in \eqref{eq:min_gen_utilization}) of 0.5 for larger thermal units, e.g., Large Simple Cycle, and Steam Turbines, and 0.05 for smaller ``peaker'' units, e.g., Small Simple Cycle.
\item We take the annualized retirement costs (\$/y) to be 10\% of the annualized FOM costs (scaled by unit capacity) so that retiring a unit results in a cost reduction of 90\% of the FOM costs.
\end{itemize}

Because the study integrates data from multiple sources with different update cycles and levels of completeness, some parameters required reconciliation or engineering assumptions. The results should therefore be interpreted as scenario-based planning insights conditional on the assembled dataset, rather than as a precise forecast of future system operation.

\section{Simulations}\label{sec:simulations}
We now present computational results for solving the CEP problems on the PR test system.  The model is formulated in \texttt{Python} using the open-source algebraic modeling language \texttt{pyomo} \cite{bynum2021pyomo,hart2011pyomo} and the open-source UC/ED modeling package EGRET \cite{knueven2020mixed}.  To decompose these large-scale optimization problems, we use a parallel implementation of Progressive Hedging (PH) with a hub-and-spoke architecture implemented in \texttt{mpi-sppy} \cite{mpi-sppy}.  This PH configuration includes incumbent finders and lower bounding methods, allowing us to compute optimality gaps that provide a certificate on the quality of our solutions (the lower the gap, the closer to provably optimal).  We decompose the 3-month time horizon by representative day; single-day subproblems are solved with the commercial optimization software Gurobi.  We run these simulations on the high performance computing (HPC) cluster \texttt{dane} at LLNL.  Each compute node has 112 2.0GHz cores and 256GB of RAM, and we assign 5 single-day scenarios per node.  

\subsection{Baseline Results without Expansion}\label{sec:noexpansion}
\begin{figure}[t]
    \centering
        \caption{Generation capacity of the 2024 baseline PR system}
    \includegraphics[width=\linewidth]{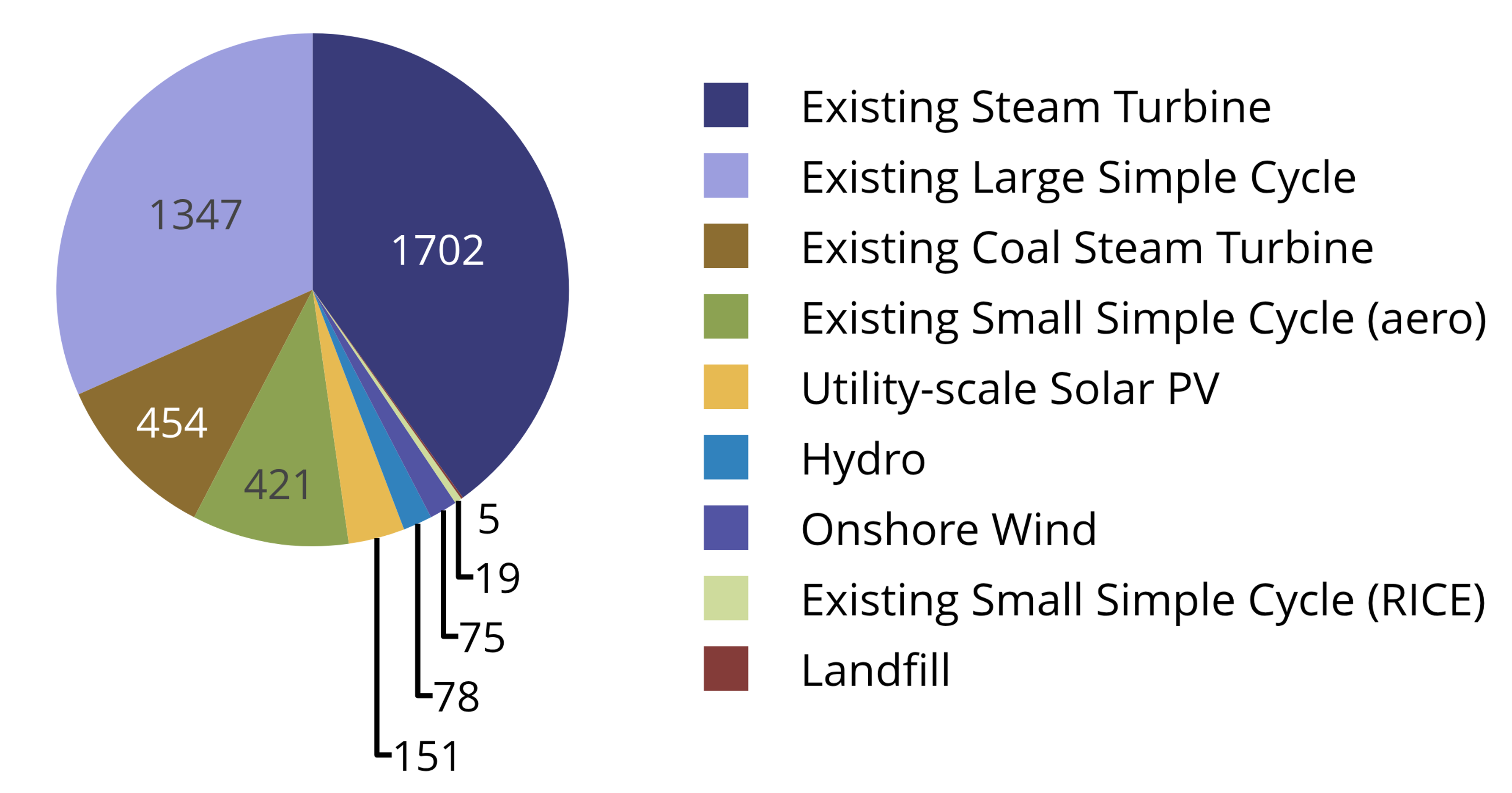}
    \label{fig:2024_piechart}
\end{figure}

\begin{figure*}[t]
    \centering
    \caption[Map of the 2024 baseline PR system]{Map of the 2024 baseline PR system. Note that some branches and generators were missing location data and are excluded from this graphic.}
\begin{center}
    \includegraphics[width=1\linewidth]{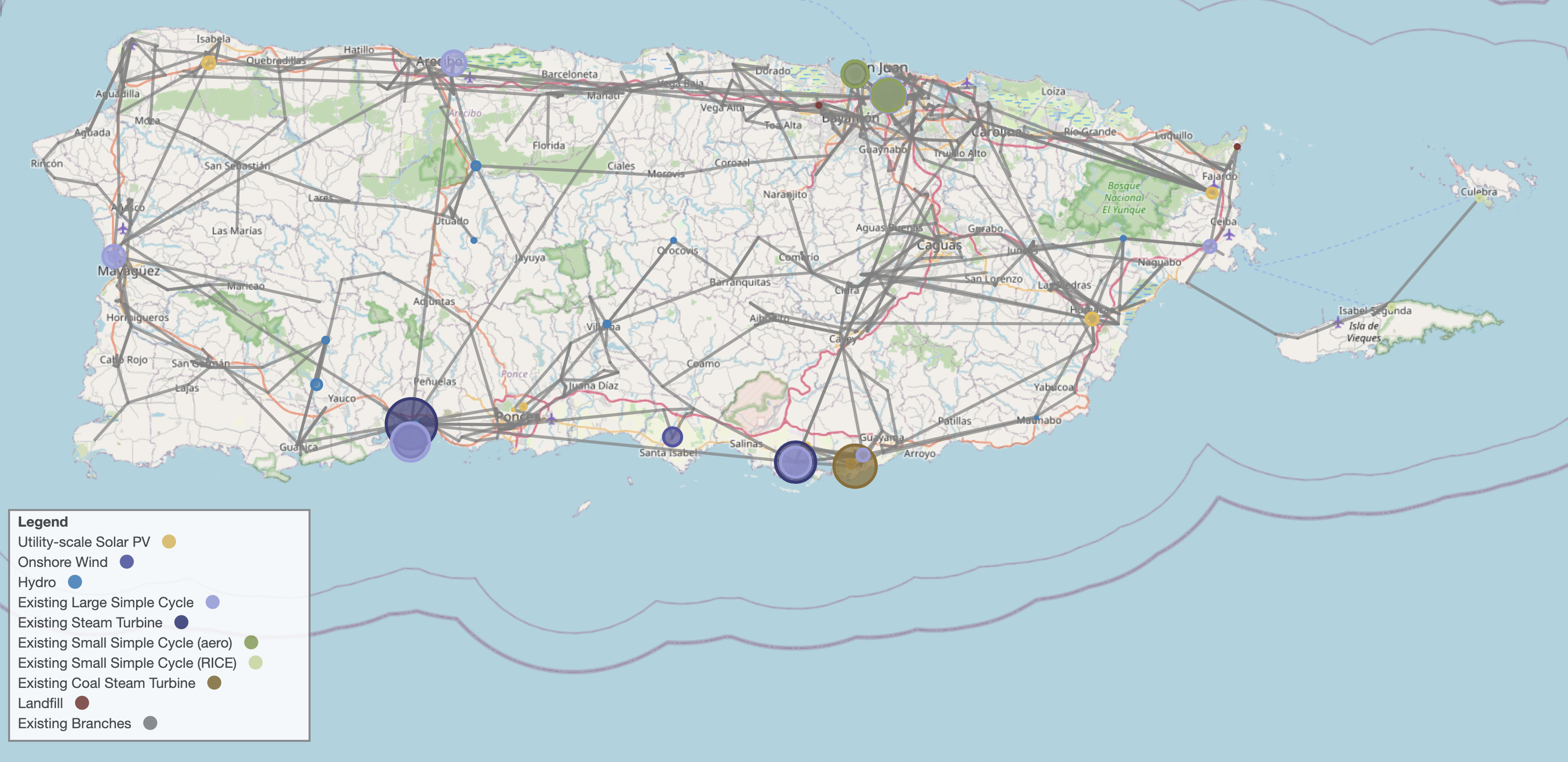}
\end{center}
\label{fig:2024_map}
\end{figure*}

First, we consider the existing PR system without allowing for new generation or storage, or retirements of existing generators.  See Figures \ref{fig:2024_piechart} and \ref{fig:2024_map} for depictions of this 2024 baseline PR system.  For these simulations, we use the 2024 time-series data for load and generator availability and 2024 fuel supply data.  Cases (E1) and (E2) use 3-month time horizons corresponding to the peak summer season, and cases (E3) and (E4) use the whole year of 2024 data. We also compare including a pipe-and-bubble representation of the transmission system to not including one, i.e.,  using a copperplate model that ignores transmission constraints.  Because we do not consider first-stage decision variables in these simulations, the problem over the entire time horizon trivially decomposes by day (each of the single-day problems is independent), and we can solve all single-day problems in less than 10 seconds to under a 1\% optimality gap.  The results from these initial simulations are given in Table \ref{tab:baseline_noexpansion}.

\begin{table*}[t]
\centering
\footnotesize
\caption{Baseline PR system results without capacity expansion}
\label{tab:baseline_noexpansion}
    \rowcolors{2}{llnlLivermoriumIce}{white}
\begin{tabular}{
>{\raggedright\arraybackslash}p{1.0cm}
>{\raggedright\arraybackslash}p{2.3cm}
>{\raggedright\arraybackslash}p{3.1cm}
>{\raggedright\arraybackslash}p{2.2cm}
>{\raggedright\arraybackslash}p{2.1cm}
>{\raggedright\arraybackslash}p{2.5cm}
>{\raggedright\arraybackslash}p{1.8cm}
}
\rowcolor{llnlElementalNavy}
\myth{Case} &
\myth{Time-series data (2024)} &
\myth{Network model} &
\myth{Production cost} &
\myth{Exp. daily load shed} &
\myth{Expected LOLH} &
\myth{Minimum reserves} \\
(E1) & May--July      & Copperplate       & 2.40 B\$/y & 521 MWh & 195 / 2208 h  & 0 MW \\
(E2) & May--July      & Pipe-and-bubble   & 2.43 B\$/y & 796 MWh & 1686 / 2208 h & 1 MW \\
(E3) & Jan--Dec.      & Copperplate       & 2.16 B\$/y & 549 MWh & 822 / 8784 h  & 0 MW \\
(E4) & Jan--Dec.      & Pipe-and-bubble   & 2.19 B\$/y & 811 MWh & 4459 / 8784 h & 1 MW \\
\end{tabular}
\end{table*}

Comparing cases (E1) and (E2), without and with the pipe-and-bubble network representation, we see that there is some congestion in the network, as observed in the increased production cost, and higher expected load shedding and loss of load hours (LOLH). We observe the same pattern when we compare cases (E3) and (E4).  When we compare the 3-month cases (E1) and (E2) to the full year cases (E3) and (E4), we see that the May to July 3-month horizon has a higher production cost, as well as a higher overall cost (including lost load and reserve shortfall penalties), than the entire year -- validating our choice to use the May to July time horizon.

We also use these existing system cases to perform verification and validation on our baseline results.  We compared our expected LOLH to LUMA's expected LOLH (as computed in their March 2025 RA study \cite{luma_ra_study_2025}) of 318 h / 2208 h for the same 3-month window (May to July 2024), which is consistent with our results -- higher than our expected LOLH for the copperplate (E1) case and lower than our expected LOLH for the pipe-and-bubble (E2) case.  For (E2), we found that generation production costs in our model were approximately 11¢/kWh, compared to an estimated 15¢/kWh based on EIA Form 861M \cite{eia861m2025apr}, where the generation portion of a PR electricity bill is estimated to be 62\% of reported residential price based on a sample LUMA bill \cite{luma_bill_charges}.  Note that our cost estimate excludes FOM costs for existing power plants when data were not available.  See the Appendix for a plot of power generation by technology type in the (E2) case.

Because these UC/ED problems without capacity expansion are much easier to solve than the full CEP problems that include first-stage investments, we also solve the entire May-July 2024 case as a single UC/ED problem, without decomposing by day.  Both of the cases -- denoted (E5) and (E6) -- solved to under 1\% optimality gaps in less than 15 minutes.  In Table \ref{tab:baseline_noexpansion_stitched}, we see that the production costs are somewhat lower for these cases, compared to their counterparts (E1) and (E2) in Table \ref{tab:baseline_noexpansion}, with less expected daily load shedding and expected LOLH.

\begin{table*}[t]
\centering
\footnotesize
\caption{Baseline PR system results, solving 3-month scenario as a single UC/ED problem}
\label{tab:baseline_noexpansion_stitched}
    \rowcolors{2}{llnlLivermoriumIce}{white}
\begin{tabular}{
>{\raggedright\arraybackslash}p{1cm}
>{\raggedright\arraybackslash}p{2.8cm}
>{\raggedright\arraybackslash}p{1.9cm}
>{\raggedright\arraybackslash}p{2cm}
>{\raggedright\arraybackslash}p{2cm}
>{\raggedright\arraybackslash}p{2.4cm}
>{\raggedright\arraybackslash}p{1.9cm}
}
\rowcolor{llnlElementalNavy}
\myth{Case} &
\myth{Network model} &
\myth{Solve time} &
\myth{Production cost} &
\myth{Exp. daily load shed} &
\myth{Expected LOLH} &
\myth{Minimum reserves} \\
(E5) & Copperplate      & 684 s  & 2.37 B\$/y & 291 MWh & 145 / 2208 h  & 0 MW \\
(E6) & Pipe-and-bubble  & 843 s & 2.41 B\$/y & 449 MWh & 1654 / 2208 h & 1 MW \\
\end{tabular}
\end{table*}

\subsection{CEP on Future PR System with \\Sensitivities to Fuel Supply and \\Baseline Fleet}\label{sec:cep_results}
Now we consider future scenarios where we allow for first-stage investment decisions, including generation and storage investments and existing thermal generator retirements.  For these simulations, we use the May to July 3-month time horizon, with generator availability data from 2024 and scaled future load corresponding to the year 2030. We consider several sensitivities for these simulations, as follows:
\begin{itemize}
    \item Three different load scenarios
    \begin{itemize}
        \item ``A'' corresponds to 2030 load case (7\% scaled \textit{increase} from 2024 load)
        \item ``B'' corresponds to elevated 2030 load case (20\% scaled \textit{increase} from 2024 load)
        \item ``C'' corresponds to reduced 2030 load case (20\% scaled \textit{decrease} from 2024 load)
    \end{itemize}
    \item Two different baseline fleet scenarios, i.e.,  what we are taking as our existing system for the CEP model
    \begin{itemize}
        \item Planned system for 2030
        \item Planned system for 2030 removing the tranches and projects funded from the Community Development Block Grant (CDBG) program
    \end{itemize}
    \item Two scenarios for future fuel supplies
    \begin{itemize}
        \item Expected ``future'' fuel supply (see Table~\ref{tab:fuel_scenarios}).
        \item Expanded ``future+'' fuel supply that also includes unlimited natural gas at Aguirre and AES
    \end{itemize}
    \item Two scenarios for what type of new generation investments we will consider
    \begin{itemize}
        \item Allow all generation types to be considered, with up to 2x the wind and solar build at existing sites.
        \item Restrict generation types to renewable-only technologies, i.e.,  wind or solar, with up to 2x wind at existing sites and unlimited solar at high-voltage buses.
    \end{itemize}
    \item Summer versus winter time-series data for the renewables-only cases
    \begin{itemize}
        \item ``Summer'' (denoted ``s'') uses May-July 2024 data, the same as all the other cases
        \item ``Winter'' (denoted ``w'') uses January-February 2024 and December 2024 data
    \end{itemize}
\end{itemize}

We additionally consider some baseline future simulations -- cases (F1-A),(F1-B), and (F1-C) -- without new investment but with retirements allowed. These (F1) simulations provide a reference point for how much generation is retired based on non-utilization in the planned future system, rather than based on replacements from more reliable new generation.

\begin{table*}[ht]
\centering
\footnotesize
\caption{Summary of future PR system sensitivities and results}
\label{tab:future_cases_all}
\rowcolors{2}{llnlLivermoriumIce}{white}
\begin{tabular}
{>{\raggedright\arraybackslash}p{1.4cm} 
 >{\raggedright\arraybackslash}p{2.0cm} 
 >{\raggedright\arraybackslash}p{1.0cm} 
 >{\raggedright\arraybackslash}p{2.0cm} 
 >{\raggedright\arraybackslash}p{0.8cm} 
 >{\raggedright\arraybackslash}p{0.8cm} 
 >{\raggedright\arraybackslash}p{0.8cm} 
 >{\raggedright\arraybackslash}p{0.8cm} 
 >{\raggedright\arraybackslash}p{0.8cm} 
 >{\raggedright\arraybackslash}p{0.8cm} 
 >{\raggedright\arraybackslash}p{0.8cm} 
 >{\raggedright\arraybackslash}p{0.8cm}} 
\rowcolor{white}
\multicolumn{4}{c}{\textbf{Problem Setup}} &
\multicolumn{8}{c}{\textbf{Results}} \\
\cmidrule(lr){1-4}\cmidrule(lr){5-12}
\rowcolor{llnlElementalNavy}
\myth{Case} &
\rothead{\myth{Baseline fleet of generation \& storage}} &
\rothead{\myth{Fuel supply}} &
\rothead{\myth{New investments allowed?}} &
\rothead{\myth{Optimality gap}} &
\rothead{\myth{Investment cost\textsuperscript{4} (M\$/y)}} &
\rothead{\myth{Production cost (B\$/y)}} &
\rothead{\myth{Overall cost\textsuperscript{5} (B\$/y)}} &
\rothead{\myth{Min. reserves (MW)}} &
\rothead{\myth{New gen (MW)}} &
\rothead{\myth{New storage (MW)}} &
\rothead{\myth{Retired gen. (MW)}} \\

(F1-A)  & 2030 planned                    & future  & None            & 1.53\%  & 0   & 2.01 & 3.49 & 493 & 0    & 0    & 2751   \\
(F1-B)  & 2030 planned                    & future  & None            & 1.99\%  & 0   & 2.31 & 4.66 & 145 & 0    & 0    &  2391   \\
(F1-C)  & 2030 planned                    & future  & None            & 1.25\%  & 0   & 1.43 & 1.99 & 726 & 0    & 0    & 3418    \\
(F2-A)  & 2030 planned                    & future  & All types       & 0.78\%  & 128 & 1.61 & 2.18 & 900 & 1689 & 0    & 3838 \\
(F2-B)  & 2030 planned                    & future  & All types       & 0.84\%  & 171 & 1.79 & 2.40 & 900 & 2239 & 0    & 3838 \\
(F2-C)  & 2030 planned                    & future  & All types       & 0.43\%  & 89  & 1.22 & 1.75 & 900 & 1157 & 0    & 3838 \\

(F3-A)  & 2030 planned                    & future+ & All types       & 0.76\%  & 130 & 1.61 & 2.18 & 900 & 1701 & 0    & 3838 \\
(F3-B)  & 2030 planned                    & future+ & All types       & 0.85\%  & 170 & 1.79 & 2.40 & 900 & 2230 & 0    & 3838 \\
(F3-C)  & 2030 planned                    & future+ & All types       & 0.44\%  & 89  & 1.22 & 1.75 & 900 & 1157 & 0    & 3838 \\

(F4-A)  & w/o tranches\textsuperscript{1} & future+ & All types       & 0.79\%  & 167 & 1.61 & 2.21 & 900 & 2211 & 0    & 3838 \\
(F4-B)  & w/o tranches\textsuperscript{1} & future+ & All types       & 0.96\%  & 212 & 1.79 & 2.44 & 900 & 2776 & 0    & 3838 \\
(F4-C)  & w/o tranches\textsuperscript{1} & future+ & All types       & 0.71\%  & 102 & 1.24 & 1.78 & 900 & 1267 & 0    & 3838 \\

(F5-As) & 2030 planned                    & future  & Renew-only\textsuperscript{3} & 12.65\% & 446 & 1.37 & 2.68 & 690 & 4586 & 32   & 3123 \\
(F5-Aw)\textsuperscript{6} & 2030 planned                    & future  & Renew-only\textsuperscript{3} & 5.22\%  & 335 & 1.23 & 2.17 & 659 & 3482 & 0    & 3493 \\
(F5-B)  & 2030 planned                    & future  & Renew-only\textsuperscript{3} & 12.80\% & 832 & 1.40 & 3.14 & 492 & 6730 & 1143 & 3144 \\
(F5-C)  & 2030 planned                    & future  & Renew-only\textsuperscript{3} & 1.88\%  & 293 & 1.03 & 1.77 & 789 & 3041 & 0    & 3743 \\

(F6-A)  &  w/o batteries\textsuperscript{2}                            &   future+       &   All types             & 1.35\%  & 176 & 1.50 & 2.12 & 825 & 2134 & 30   & 3838 \\
(F6-B)  &  w/o batteries\textsuperscript{2}                             &   future+       &   All types                & 2.66\%  & 200 & 1.71 & 2.36 & 729 & 2308 & 121  & 3763 \\
(F6-C)  &  w/o batteries\textsuperscript{2}                              &   future+       &   All types               & 0.81\%  & 92  & 1.14 & 1.68 & 868 & 1182 & 1    & 3838 \\
\end{tabular}
\vspace{0.5em}
\raggedright

\textsuperscript{1}w/o tranches = 2030 planned w/o tranches, CDBG\\
\textsuperscript{2}w/o batteries = 2030 planned w/o batteries\\
\textsuperscript{3}Renew-only = only allow new investments in wind, solar, and BESS\\
\textsuperscript{4}This investment cost includes only the cost of new generation and storage investments, not the cost of thermal generator retirements.\\
\textsuperscript{5}The overall cost is from Objective (\ref{eqn:objective}) and includes load shed penalty (\$30{,}000/MWh), reserve shortfall penalty (\$2{,}000/MWh), generator under-utilization penalty (\$2{,}000/MWh), and retirement costs, in addition to investment and production costs.\\
\textsuperscript{6}The (F5-Aw) case uses winter time series data and is not directly comparable to the other cases.
\end{table*}

\begin{figure*}
    \centering
        \caption[Example of first-stage decisions from the CEP model]{Example of first-stage decisions from the CEP model.  For this F2-A case, the baseline system is the 2030 planned system.  The CEP model decides which thermal generator retirements and new generation and storage investments will yield the most optimal final system composition, in terms of cost and reliability, shown on the right.}
    \includegraphics[width=0.9\linewidth]{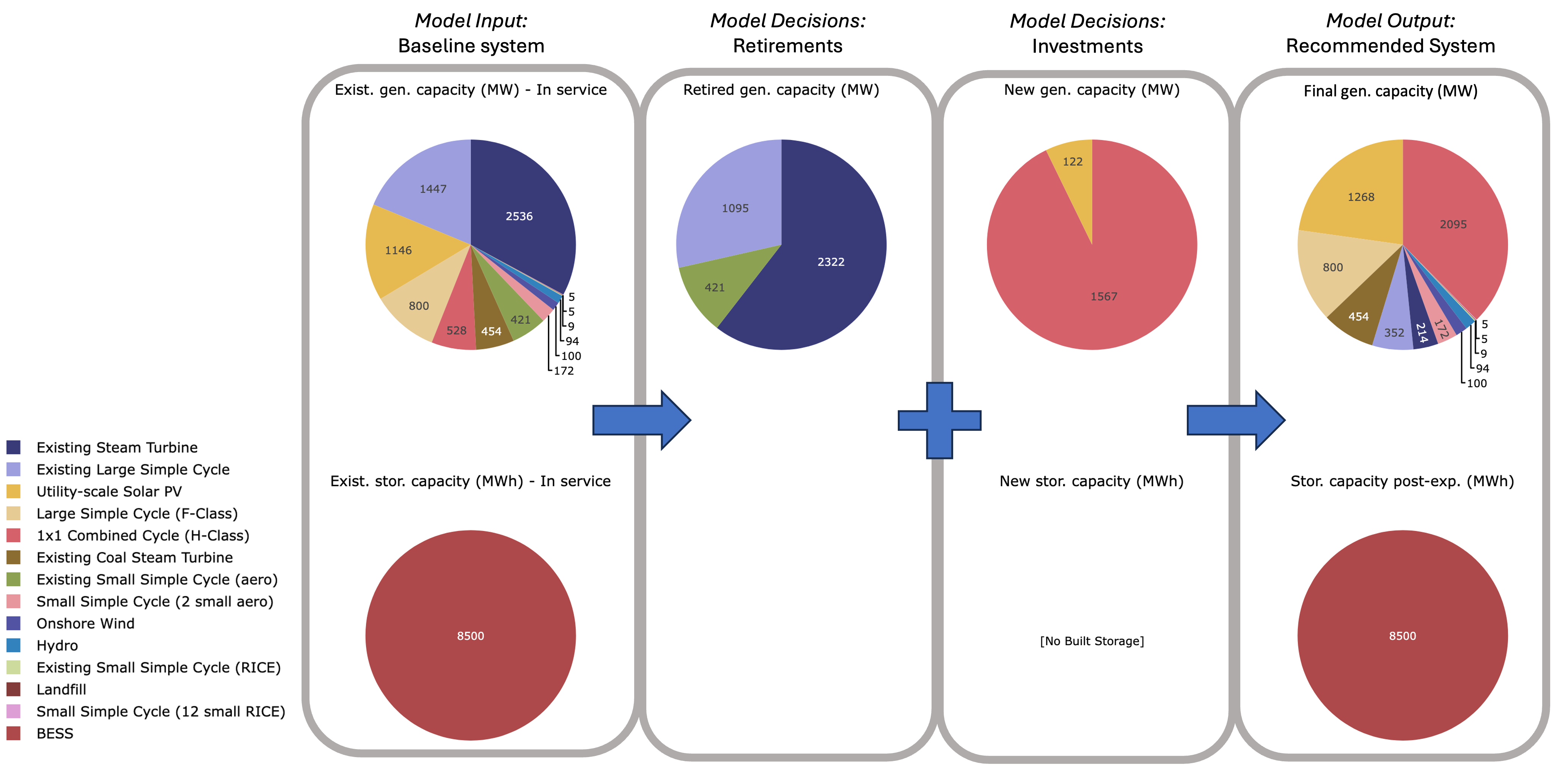}
    \label{fig:F2A}
\end{figure*}

\begin{figure}[!t]
    \centering
        \caption[Sensitivities of first-stage investment decisions to different baseline fleets and fuel supplies]{Sensitivities of first-stage decisions to different baseline fleets and fuel supplies, i.e., cases (F2-A) to (F4-A). For these cases, we saw that the model was more sensitive to baseline fleet than to fuel supply.  Without the existence of the tranches or CDBG projects, the model chooses to build an additional H-class CC unit.}
    \includegraphics[width=0.6\linewidth]{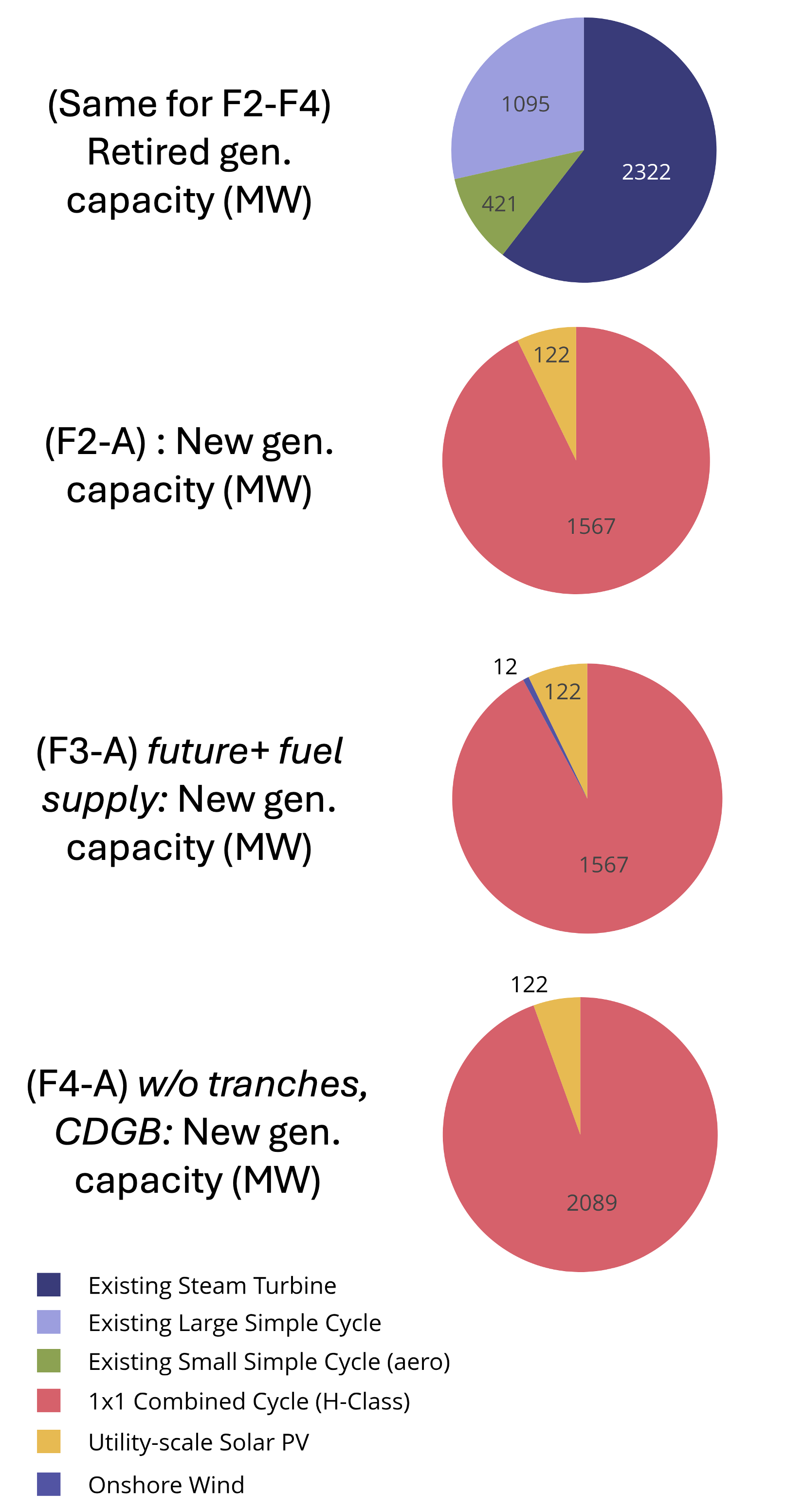}
    \label{fig:F2-F4A_sensitivities}
\end{figure}

\begin{figure}[!t]
    \centering
        \caption[Sensitivities of first-stage investment decisions to future load scenarios A, B, and C]{Sensitivities of first-stage decisions to future load scenarios A, B, and C.  For these (F2) cases, when we increase load from A to B, we see the addition of 1 more H-class CC unit.  When we decrease load from A to C, we see the reduction of one H-class unit.}
    \includegraphics[width=0.6\linewidth]{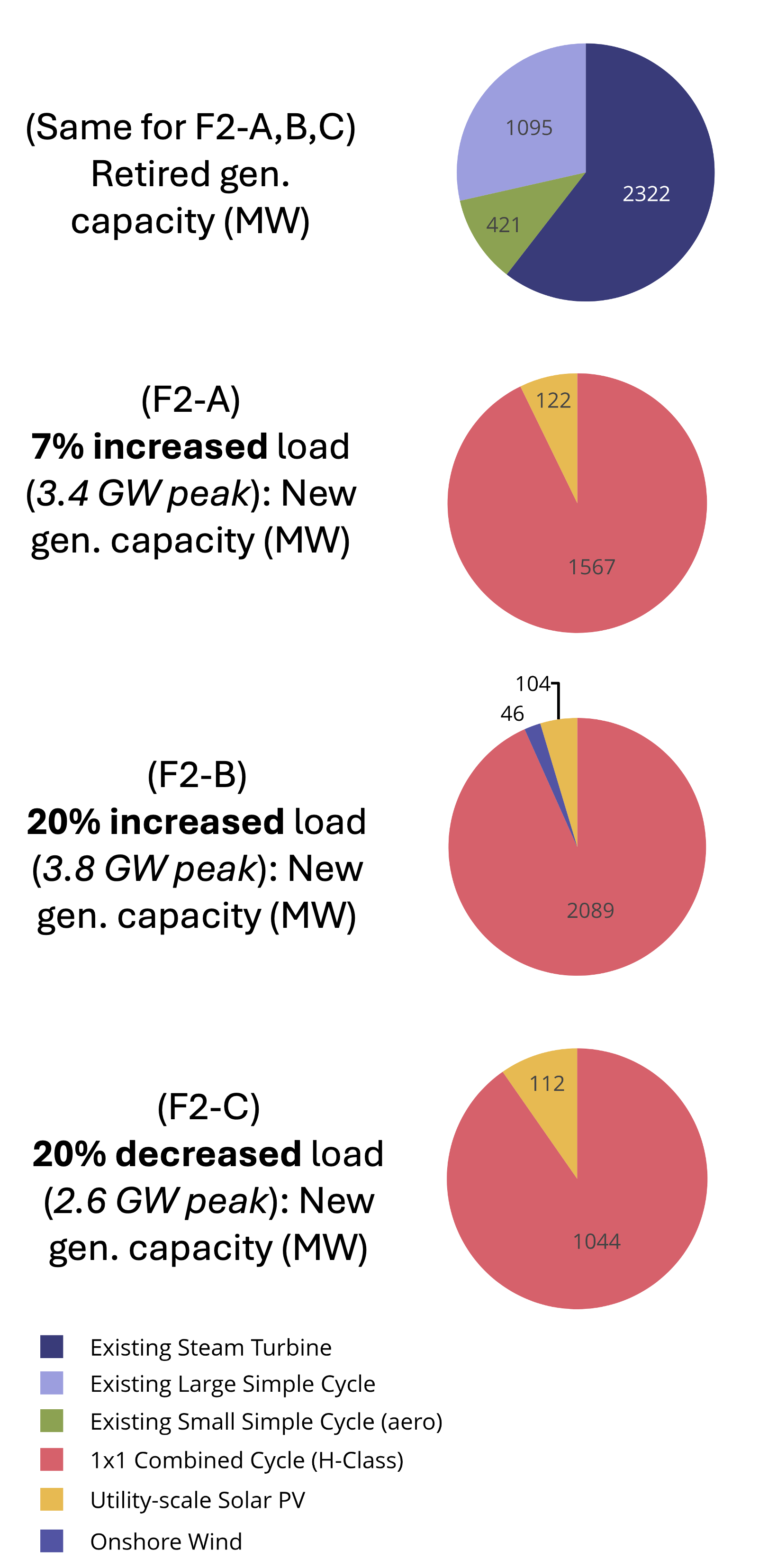}
    \label{fig:load_sensitivity}
\end{figure}

The results from these different sensitivities are summarized in Table \ref{tab:future_cases_all}. Some detailed results from Table \ref{tab:future_cases_all} are visualized in Figures \ref{fig:F2A} through \ref{fig:load_sensitivity}.  Note that we set a 20-hour wall clock time limit for PH in our CEP simulations. However, cases (F2) through (F4) all solved within 10 minutes to optimality gaps below 1\%.  For the renewables-only (F5) cases, we saw larger optimality gaps of up to 13\%.  Thus, these renewables-only results should be interpreted with some caution, as they could be somewhat far from the optimal solution -- although we are unsure if there is an upper or lower bounding issue (or both).

Relative to the modeled present-day (2024) system (E2), the modeled future (2030) system (F1-A) with only planned generation has substantially improved bulk-system performance, reducing production costs by 17\% and eliminating load shedding in the scenarios considered.  This suggests that planned repairs and additions should materially improve system reliability. However, these planned improvements do not ensure an adequate reserve margin after accounting for the retirement of non-utilized units. 


For all future cases, including the no-expansion (F1) cases, we saw no load shedding.  In all the CEP cases where we allowed new thermal generation to be built, i.e., (F2)-(F4), we met our reserve margin target of 900 MW for all time periods.  In the renewables-only cases (F5), we saw inadequate reserves as low as 492 MW in (F5-B). However, these cases still had much more reserves compared to the present-day system with 0-1 MW of reserves (see Table \ref{tab:baseline_noexpansion}). Comparing (F2-A) to (F1-A), we can see that when new investments are allowed in (F2-A), we observe approximately 1.1 GW of additional generation retirement, indicating that some of the thermal additions serve as replacements for unreliable and high-cost legacy units.  Additionally, all the CEP cases with new thermal generation allowed -- cases (F2)-(F4) -- saw the same quantity of retired generation, irrespective of load scenario, implying that some set of unreliable thermal generators should be replaced with updated technology for reliability rather than load growth concerns.

In the renewables-only cases (F5), we saw a reduced amount of thermal generator retirement, as expected.  
Under the study assumptions, restricting new builds to renewables and storage increased total modeled system cost relative to cases that allowed new thermal construction.
On average across the three summer load scenarios, the overall costs were about 20\% higher for the renewables-only cases (F5) compared to the thermals-allowed case (F2).

 Looking at the different sensitivities, we find that the fuel expansion from case (F2) to case (F3) has almost no impact on the results. The fleet change from (F3) to (F4), with the removal of tranches and CDBG projects, results in somewhat greater costs and higher system investment, most apparent in the (F4-A) and (F4-B) load cases.  We see that costs increase and more generation is built in the higher load cases ``B'' and vice versa for the reduced load cases ``C,'' as expected.  For example, the increased load from case (F2-A) to case (F2-B) results in the addition of one more H-class combined cycle unit, as seen in Figure \ref{fig:load_sensitivity}.  In the renewables-only summer vs. winter sensitivity, we see that the costs and investment quantities were larger for the summer case.  Our finding that the summer renewables-only case is the more ``stressed'' case is consistent with this observation in our data: while solar PV is 20\% more available in the summer than the winter, the loads are 23\% higher.

    \begin{figure}[!t]
    \centering
        \caption[Sensitivities of first-stage investment decisions to allowable technology types for new generation (thermals allowed vs. renewables only)]{Sensitivities of first-stage decisions to allowable generation types (thermal allowed vs. renewables only).  For these cases, we saw about 700 MW less thermal generator retirement in the renewables-only case (F5-As) compared to the thermals-allowed case (F2-A).  We see that (F5-As) requires over 4 GW of new nameplate capacity while (F2-A) requires only about 1.7 GW.}
    \includegraphics[width=1\linewidth]{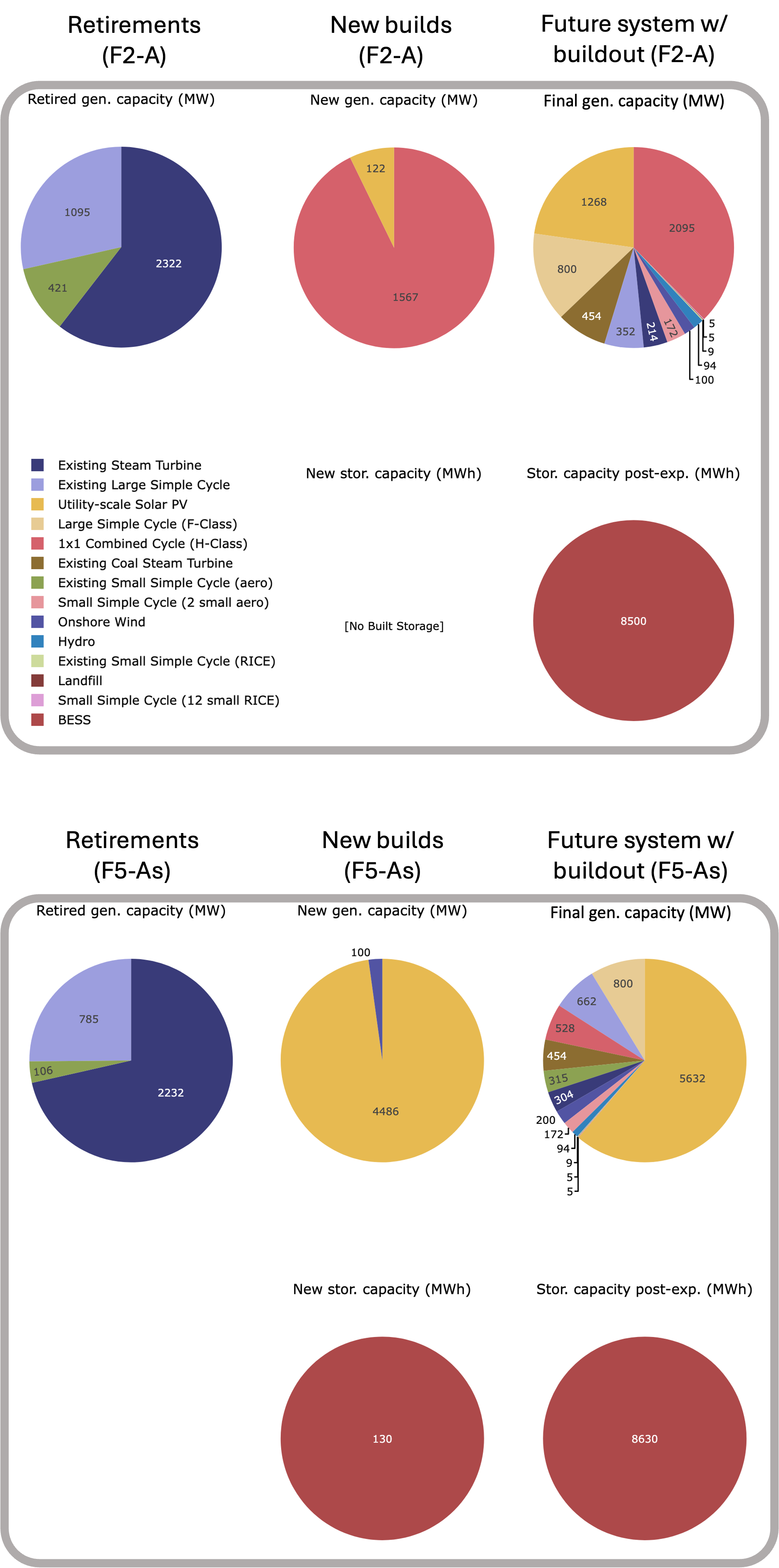}
    \label{fig:F2vsF5}
\end{figure}

\subsection{CEP on Future PR System with \\Sensitivities to Battery Fleet}
We also considered removing all the batteries from the baseline PR system, specifically all 2.1 GW of planned BESS projects from the 2030 system (note that there are no utility-scale batteries in our 2024 system).  We denote these cases with an (F6) prefix and note that they are equivalent to the (F3) cases without the planned BESS in the baseline fleet.  These results are shown in Table \ref{tab:future_cases_all}.  All these cases solved to under a 3\% optimality gap within the 20-hour time limit.
  
When we compare the planned versus no-battery baseline fleets for all three load cases, we see that there are somewhat higher amounts of generation investment and small amounts of storage investment for the no-battery (F6) cases compared to the planned (F3) cases.  We also find that the (F3) cases have higher overall costs compared to the (F6) cases. This is due to the fact that FOM costs associated to these planned storage projects are eliminated, resulting in a cost savings of 0.87 B\$/y. This analysis does not capture the investment costs for these planned BESS units, and our CEP model does not build much, if any, storage in the scenarios where new thermal units can be built (see above section). 
Overall, our results suggest the full planned BESS fleet is not required to minimize cost when new thermal capacity is available. This result should be interpreted cautiously because the model does not fully capture dynamic stability, contingency response, or other ancillary service value streams that can be provided by storage.


\subsection{CEP on Future PR System with No Existing or New Thermal Units}
\begin{figure}[!t]
    \centering
        \caption[Sensitivities of first-stage decisions to the removal of all existing thermal units]{Sensitivities of first-stage decisions to the removal of all existing thermal units.  For these cases, we saw about 4x the solar PV buildout and over 200x the BESS buildout in the (F7) case that has no thermal units compared to the (F5) case that preserves existing thermal units.  Neither of these cases allows for new thermal units to be built.}
    \includegraphics[width=\linewidth]{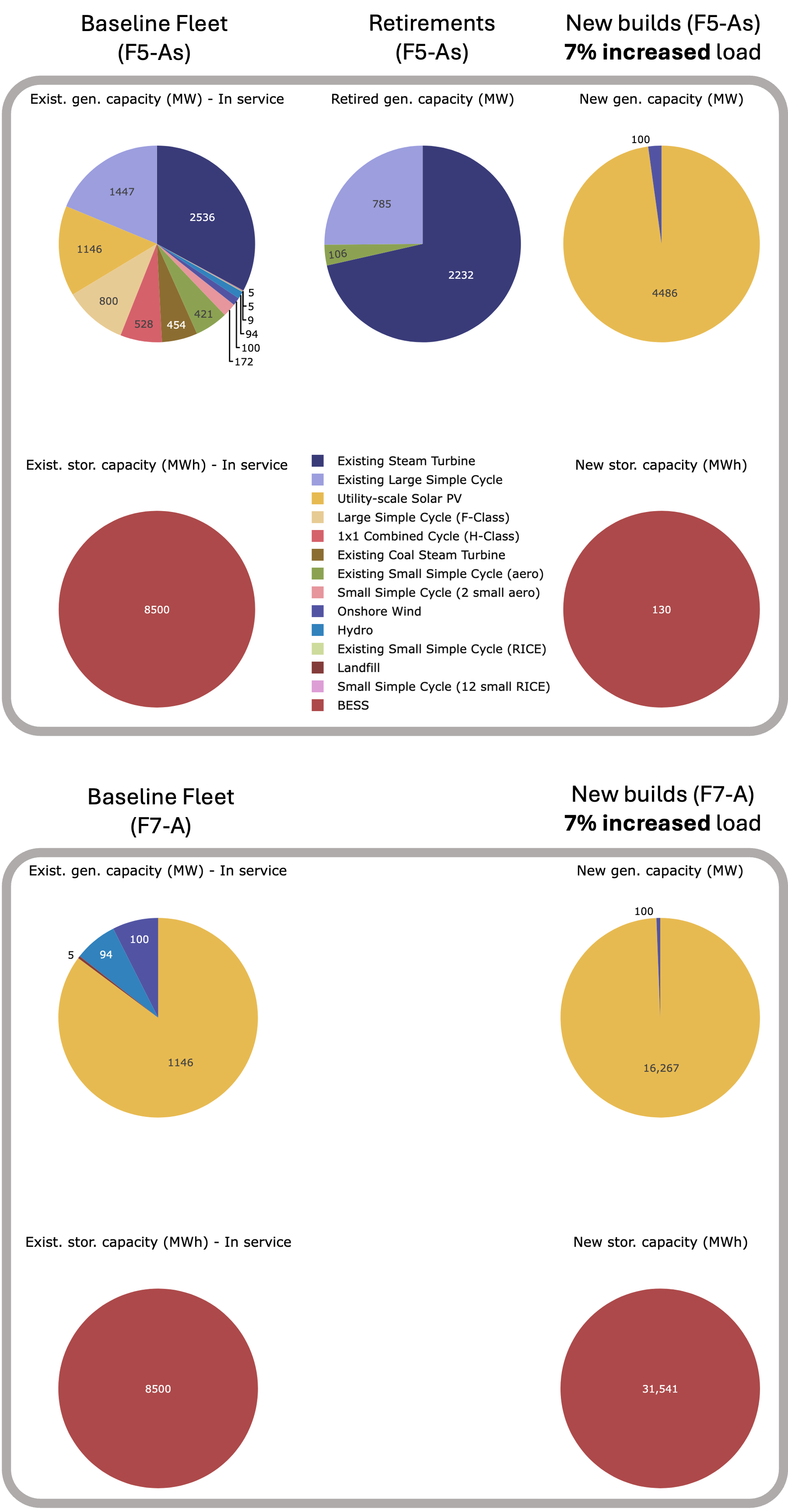}
    \label{fig:nothermals}
\end{figure}

Next, we consider the scenario where all existing thermal units are removed from the baseline system and no new thermal units are allowed to be built, in order to model a policy of 100\% renewable-only generation.  We label these cases (F7) and note that they are equivalent to the (F5) cases without the existing thermal generators in the baseline fleet.  The results for this sensitivity are given in Table \ref{tab:nothermals_results}. All these cases solved to within a 3\% optimality gap within the 20-hour time limit.  These results show that investment costs increase substantially and load shedding remains significant 
if we have no existing or new thermal units in the future system.  For example, comparing cases (F3-A) and (F7-As), we see that the sum of investment and production costs are 71\% higher when we have no thermal units in the future system.  When we look at overall costs that include penalties for load shedding and reserve shortfall, we see that the cost is 35 times higher for the no-thermals case (F7-As) compared to (F3-A).  We see much higher solar PV and BESS build-out in these (F7) cases compared to the (F5) cases, as shown in Figure \ref{fig:nothermals} for load case A.  This sensitivity is the most directly comparable case to the 100\% renewables futures considered in PR100. However, we note that PR100 included some dispatchable renewable energy options such as biomass and multiple extended-duration (up to 10 hours) storage technologies that our model does not consider.

\begin{table*}[t]
\centering
\footnotesize
\caption{Summary of future PR system results with no-thermal sensitivities}
\label{tab:nothermals_results}
    \rowcolors{2}{llnlLivermoriumIce}{white}
\begin{tabular}{
>{\raggedright\arraybackslash}p{1.1cm}
>{\raggedright\arraybackslash}p{1.2cm}
>{\raggedright\arraybackslash}p{1.8cm}
>{\raggedright\arraybackslash}p{1.8cm}
>{\raggedright\arraybackslash}p{1.8cm}
>{\raggedright\arraybackslash}p{1.8cm}
>{\raggedright\arraybackslash}p{1.0cm}
>{\raggedright\arraybackslash}p{1.5cm}
>{\raggedright\arraybackslash}p{1.3cm}
>{\raggedright\arraybackslash}p{1.3cm}
}
\rowcolor{llnlElementalNavy}
\myth{Case} &
\myth{Optimal-ity gap} &
\myth{Investment cost (M\$/y)} &
\myth{Production cost (B\$/y)} &
\myth{Overall cost (B\$/y)} &
\myth{Avg. daily load shed} &
\myth{LOLH out of 2160 h} &
\myth{Min. reserves (MW)} &
\myth{New gen. (MW)} &
\myth{New storage (MW)} \\
(F7-As) & 0.99\% & 2850 & 0.123 & 75.54  & 5.2 GWh  & 802 h  & 0 & 16367 & 7885 \\
(F7-Aw) & 1.24\% & 2750 & 0.123 & 69.16  & 4.6 GWh  & 746 h  & 0 & 16448 & 7209 \\
(F7-B)  & 0.81\% & 3000 & 0.124 & 132.51 & 10.4 GWh & 1253 h & 0 & 16397 & 8749 \\
(F7-C)  & 2.17\% & 2450 & 0.123 & 27.20  & 809 MWh  & 138 h  & 0 & 16313 & 5444 \\
\end{tabular}
\end{table*}

\subsection{CEP on Near-term PR System with Sensitivity in Types of Generation Available for New Construction}\label{sec:gradual_buildout}
Finally, we consider the optimal buildouts for the years before 2030.  Based on the observation that it is unlikely for a new CC unit to come online before 2031 due to long lead times, we examine how the system buildout would differ if we took different years between 2026 and 2031 as the baseline system.  Note that this is not a multi-stage approach in that we do not optimize over all the years gradually-- rather, we are taking a one-shot approach at each year.  

\begin{figure}[!t]
\centering
\caption[Sensitivities of the first-stage investment decisions to future buildout year, from 2027 to 2031]{Sensitivities of the first-stage decisions to future buildout year, from 2027 to 2031.  While CC (H-class) is part of the final optimal portfolio in 2031, it is not available for nearer-term buildouts in these scenarios.  A mix of different technologies, including RICE, solar and BESS, may be more suitable for near-term system targets.}
    \includegraphics[width=0.8\linewidth]{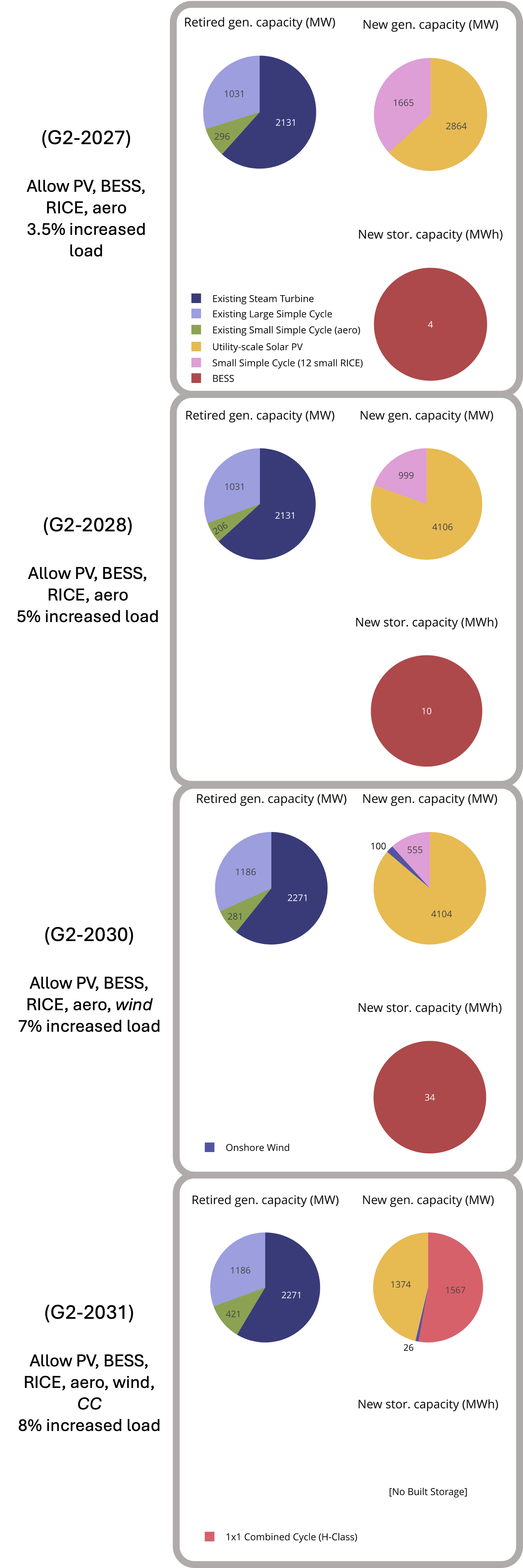}
    \label{fig:incremental_study}
\end{figure}

For these scenarios, we model load growth from 2024 to 2031, interpolating load linearly between the years based on 7\% increased load in 2030 compared to the load in 2024.  The baseline system between the years changes as different planned units come online, as reflected by their in-service dates in the data.  Note that we are considering units online for a given year if they are online before August 1 of that year.  We capture different lead times for different types of technologies: we allow RICE, solar, aero, and BESS to come online in 2027, wind in 2030, and CC in 2031.  For all these cases, we use the ``future'' fuel supply scenario.  In a difference from other cases (except for the renewables-only cases), we allow unlimited solar at buses that are considered for expansion.  The results are given in Table \ref{tab:incremental_study} and visualized in Figure \ref{fig:incremental_study}.  We see in these results that system costs decrease as planned units come online and more types of generation are allowed to be built.

\begin{table*}[t]
\centering
\footnotesize
\caption{Summary of near-term PR system sensitivities}
\label{tab:incremental_study}
    \rowcolors{2}{llnlLivermoriumIce}{white}
\begin{tabular}{
>{\raggedright\arraybackslash}p{1.6cm}
>{\raggedright\arraybackslash}p{1.6cm}
>{\raggedright\arraybackslash}p{1.6cm}
>{\raggedright\arraybackslash}p{1.6cm}
>{\raggedright\arraybackslash}p{1.6cm}
>{\raggedright\arraybackslash}p{1.4cm}
>{\raggedright\arraybackslash}p{1.2cm}
>{\raggedright\arraybackslash}p{1.2cm}
>{\raggedright\arraybackslash}p{1.2cm}
}
\rowcolor{llnlElementalNavy}
\myth{Case} &
\myth{Investment cost (M\$/y)} &
\myth{Production cost (B\$/y)} &
\myth{Overall cost (B\$/y)} &
\myth{Avg. daily load shed (MWh)} &
\myth{Min. reserves (MW)} &
\myth{New gen. (MW)} &
\myth{New storage (MW)} &
\myth{Retired gen. (MW)} \\
(G2-2026) &    0  & 2.20 & 4.75 & 8.5  &   0 &    0   & 0 &    0 \\
(G2-2027) &  500  & 1.33 & 2.36 & 0    & 726 & 4529 & 1 & 3458 \\
(G2-2028) &  530  & 1.31 & 2.40 & 0    & 715 & 5105 & 3 & 3368 \\
(G2-2030) &  480  & 1.30 & 2.30 & 0    & 665 & 4759 & 8 & 3738 \\
(G2-2031) &  250  & 1.37 & 2.06 & 0    & 900 & 2966 & 0 & 3878 \\
\end{tabular}
\end{table*}

\section{Conclusions}
\label{sec:conclusions}

This paper presents a CEP model and data workflow for analysis of Puerto Rico's power grid.  The corresponding software tool that was developed as part of this work was used to analyze expansion decisions across a range of contexts and can continue to be used as a modeling support tool for decision makers looking to improve Puerto Rico's power grid infrastructure.
This tool and resulting analysis differs from previous capacity expansion work on the PR system in three important ways:
(1) Our CEP model is nodal with high spatial and temporal resolution; (2) it includes high-resolution modeling of the operational details of thermal generators, including unit commitment and regional constraints on fuel supply; and (3) it captures weather variance through numerous representative day samples.

Across the planning scenarios examined, least-cost portfolios that maintain modeled bulk-system reliability consistently include roughly 1.5 GW or more of new H-class combined cycle capacity, in addition to planned projects.
In the model, these units primarily replace unreliable existing thermal capacity.
Preliminary results suggest that this CC investment restores a robust reserve margin and enables consumer demand to be met with low-cost, flexible generation.

A direct assessment of the impact of an RPS policy was not part of the study. However, 
in this study, we found that restricting new builds to renewable generation and storage increased the combined investment and operational costs up to 14\%\footnote{
This 14\% cost increase corresponds to the 20\% increased load future scenario; penalty-based costs like reserve shortfall are not included in these numbers.
} relative to cases that allowed new thermal capacity.
In the case where the future fleet is restricted to renewable-only sources, the cost increase was much higher at 71\% for the projected future load scenario.  These numbers may increase if additional BESS capacity is found to be necessary to provide system stability support upon further analysis.  

Our results also suggest that the planned BESS deployment may not be cost-minimizing when new thermal construction is available.  However, it is important to recall that this study may not capture all the value provided by BESS since stability and contingency analysis were not conducted.

These findings should not be interpreted as a single recommended build plan. Rather, they identify recurring features of least-cost portfolios across the scenarios and assumptions evaluated in this study. As implied above, further studies are needed for a more definitive assessment. We discuss recommended next steps below.

\subsection{Suggestions for Future Work}

Although it is expected that the inclusion of stochastic outage scenarios and a more accurate representation of the operational characteristics of thermal power plants will lead to an investment plan that is better aligned with the reliability needs of the island, this work does not eliminate the need for a resource adequacy (RA) study.
An RA study should be performed to refine the results of this analysis and assess the need for complementary investments.

Our work demonstrates how different variable operational scenarios can be handled by this  CEP tool. As a step that could be subsequent or parallel to the RA study, a richer set of operational scenarios could be identified to better represent the most relevant operational conditions of the system. This may include a broader time horizon, including contingencies, extreme weather conditions, more forced outages, etc.

A valuable improvement to this work would be the inclusion of more realistic power flow constraints that better capture the physical limitations of the transmission system. While this direction was explored during this project, we found that it was not scalable to solve the CEP problem on our high-resolution system with even linearized power flow constraints.  Additionally, we believe some improvement of the representation of spinning reserves and ancillary services could be made to enhance the operational modeling capabilities of the tool without significantly complicating the model.
Increasing modeling detail must be undertaken with care to not compromise the computational tractability of the model, but ongoing research aims to make large-scale CEP problems scalable with high levels of operational realism.

Finally, this model did not consider the uncertainty around the ability to deploy the recommended resources.
Different technologies or sizes of power plants may face different social acceptance and supply chain difficulties.
Introducing a multi-stage planning model that captures such uncertainties would be of interest to decision makers. This would likely necessitate simplifying other parts of the model to make computational tractability possible but would be a path worth pursuing.

\section*{AI Disclaimer}
Artificial Intelligence was only used for minor editing (for spelling and grammar but not content) and to format the tables. The paper has been verified by the authors, who are solely responsible for its correctness.

\section*{Acknowledgements}
The authors would like to thank our program manager at DOE's Office of Electricity (OE), Jose Benitez Torres, for many helpful discussions and valuable feedback over the numerous iterations of this work.  We would also like to acknowledge Alex Nassif, Marcelo Elizondo, Alok Kumar Bharati, and Jason Fuller at PNNL, who have run power flow and contingency analysis of our model outputs and provided constructive feedback on our work.  Finally, we would like to acknowledge DOE and PREPA staff, respectively Jeffrey Hoffmann and Ezequiel Nieves Ayala in particular, who provided valuable data inputs to our model. The authors gratefully acknowledge Gurobi Optimization for providing a free academic license for their software, which was used in this work. 

\bibliographystyle{IEEEtran}
\bibliography{refs}

\section{Appendix}
\subsection{Supplementary Figures and Tables}
\begin{figure}[!t]
    \centering
        \caption[Comparison of capacity factor data for different solar PV and wind plants in Puerto Rico's power system]{Comparison of capacity factor data for different solar PV and wind plants in Puerto Rico's power system.  The numbers from the PLEXOS PCM were sometimes much higher than the numbers provided by LUMA corresponding to actual observed generation.  We used the data from LUMA (``actual'') for our analysis.}
    \includegraphics[width=0.9\linewidth]{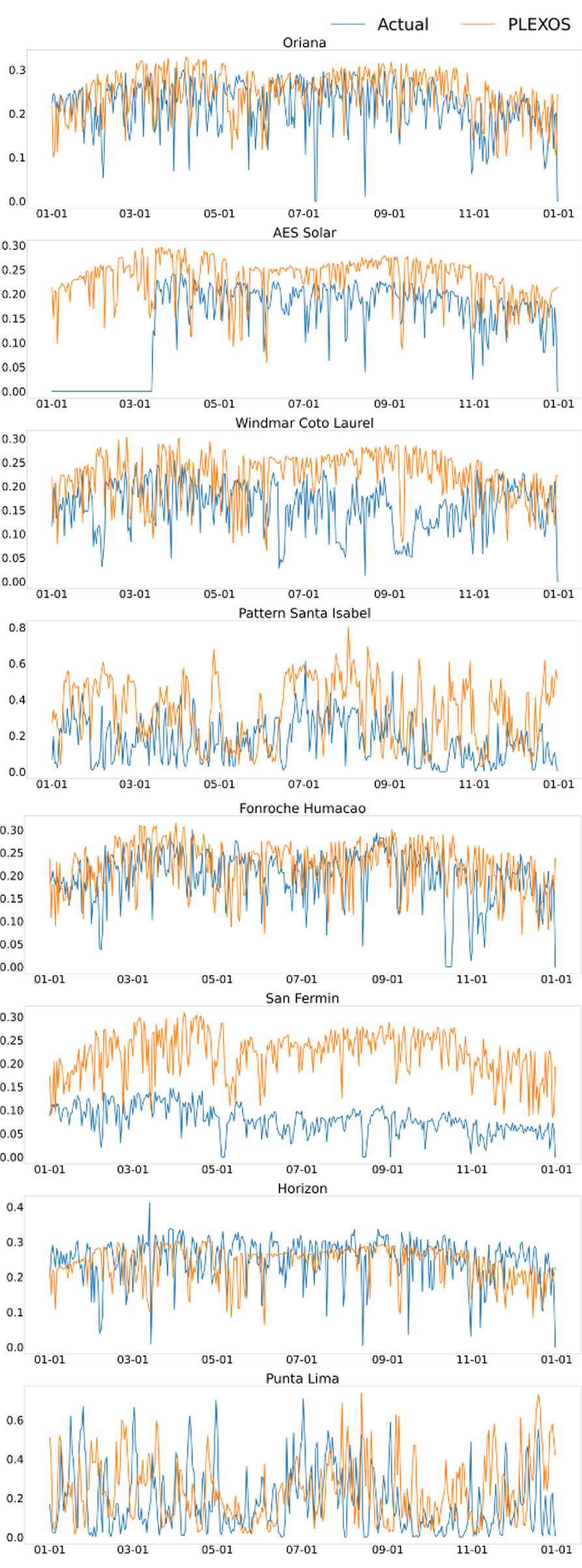}
    \label{fig:gen_avail}
\end{figure}

\begin{figure*}
    \centering
        \caption{One day power stackgraph for existing system in the (E2) test case.  We can see that the power generated is primarily from steam turbines, followed by large simple cycle.}
    \includegraphics[width=\linewidth]{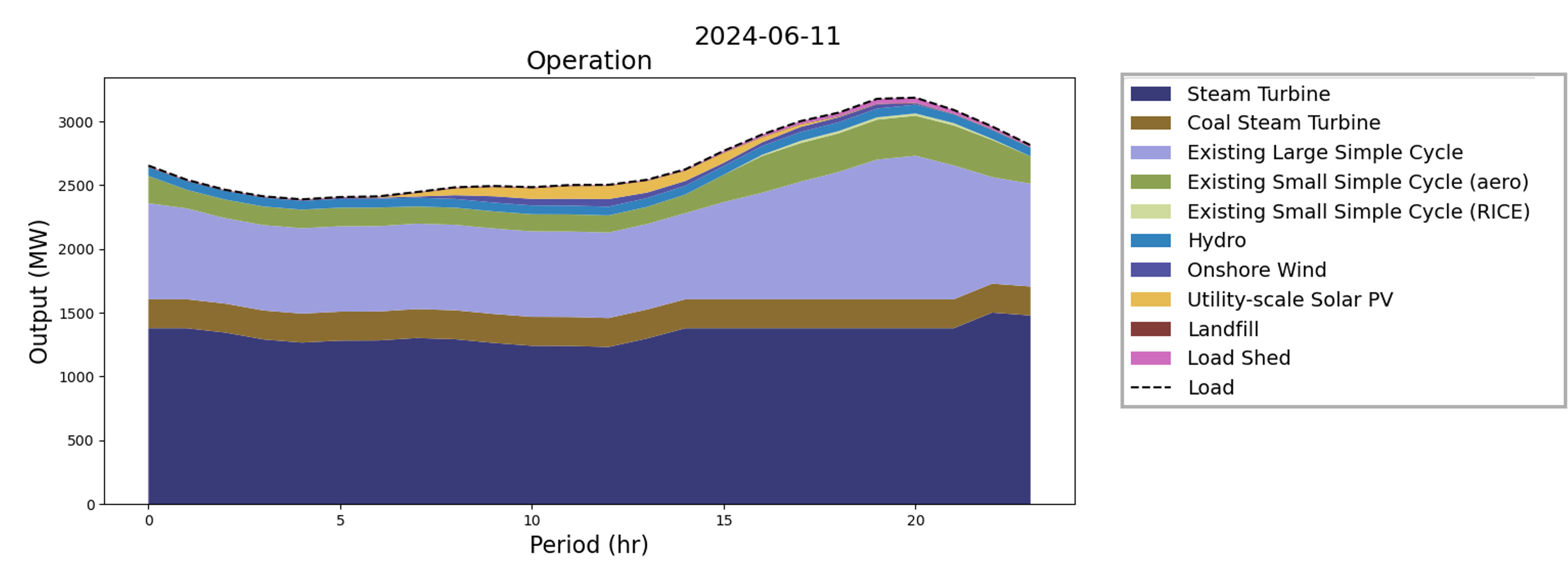}
    \label{fig:E2_stack}
\end{figure*}

See Table \ref{tab:acronyms} for a list of the acronyms and abbreviations used in this report.  See Table \ref{tab:costs_new_gen} for the cost assumptions used for new generators. See Figure \ref{fig:gen_avail} for a comparison of capacity factor data from the PLEXOS PCM and LUMA actual observed generation. See Figure \ref{fig:E2_stack} for a plot of existing system operation on the highest load day (June 11, 2024 with 3.1 GW) in the (E2) test case. 

\begin{table}[!t]
\centering
\footnotesize
\caption{Acronyms and Abbreviations}\label{tab:acronyms}
\rowcolors{2}{llnlLivermoriumIce}{white}
\begin{tabular}{>{\raggedright\arraybackslash}p{1.4cm} >{\raggedright\arraybackslash}p{6.4cm}}
\rowcolor{llnlElementalNavy}
\myth{Acronym} & \myth{Meaning} \\
AES  & AES Puerto Rico (coal-fired power plant) \\
AEO  & Annual Energy Outlook \\
BESS & Battery Energy Storage System \\
 BLS                    & Bureau of Labor Statistics \\
 CC                     & Combined Cycle (generation technology) \\
  CCGT                   & Combined Cycle Gas Turbine \\
 CDBG                   & Community Development Block Grant \\
 CEP                    & Capacity Expansion Planning \\
 CEPCI                  & Chemical Engineering Plant Cost Index (capex CPI multipliers) \\
  CS                     & Costa Sur\\
 DER                    & Distributed Energy Resource \\
 DOE                    & U.S. Department of Energy \\
  ED                     & Economic Dispatch \\
 EIA                    & U.S. Energy Information Administration \\
 EGRET                  & Electricity Grid Real-time and Expansion Toolkit (Pyomo-based UC/ED package) \\
  F-class                & “F-class” gas turbine technology (frame type) \\
 FE                     & Office of Fossil Energy (DOE-FE) \\
 FOM                    & Fixed Operations and Maintenance \\
 H-class                & “H-class” gas turbine technology (frame type) \\
  HH                     & Henry Hub (reference price for natural gas) \\
 HPC                    & High Performance Computing \\
 HR                     & Heat Rate \\
 HV                     & High Voltage \\
 IRP                    & Integrated Resource Plan\\
  LLNL                   & Lawrence Livermore National Laboratory \\
   LNG                    & Liquefied Natural Gas \\
 LOLH                   & Loss of Load Hours \\
 LUMA                   & LUMA Energy (Puerto Rico T\&D system operator) \\
 MILP                   & Mixed-Integer Linear Program \\
 NEPR                   & Puerto Rico Energy Bureau (Negociado de Energ\'ia de Puerto Rico) \\
 NETL                   & National Energy Technology Laboratory \\
 NG                     & Natural Gas \\
 OE                     & DOE Office of Electricity \\
   P3A & Puerto Rico Public-Private Partnerships Authority \\
 PCM                    & Production Cost Model \\
 PH                     & Progressive Hedging (stochastic optimization algorithm) \\
 PI                     & AVEVA PI Vision dashboard\\
 PLEXOS                 & Commercial energy modeling software \\
 PNNL                   & Pacific Northwest National Laboratory \\
 PR                     & Puerto Rico \\
 PR100                  & DOE ``Puerto Rico 100\%'' Study \\
 PREPA                  & Puerto Rico Electric Power Authority \\
 PS                     & Palo Seco \\
 PSS/E                  & Power System Simulator for Engineering (Siemens PTI software) \\
 PV                     & Photovoltaic \\
 RA                     & Resource Adequacy \\
 RAW                    & PSS/E RAW file format (power flow data) \\
 RICE                   & Reciprocating Internal Combustion Engine \\
 RPS                    & Renewable Portfolio Standard \\
 SJ                     & San Juan\\
 T\&D & Transmission and Distribution \\
 UC                     & Unit Commitment \\
 USD                    & United States Dollar \\
 VOLL                   & Value of Lost Load \\
 VOM & Variable Operations and Maintenance \\
\end{tabular}
\end{table}

\begin{table*}
\centering
\footnotesize
\caption{Cost assumptions for new generation and storage}
\label{tab:costs_new_gen}
    \rowcolors{2}{llnlLivermoriumIce}{white}
\begin{tabular}{
  >{\raggedright\arraybackslash}p{2.6cm}
  >{\raggedright\arraybackslash}p{1.4cm}
  >{\raggedright\arraybackslash}p{1.4cm}
  >{\raggedright\arraybackslash}p{2.2cm}
  >{\raggedright\arraybackslash}p{2cm}
  >{\raggedright\arraybackslash}p{2.0cm}
  >{\raggedright\arraybackslash}p{3.0cm}
}
\rowcolor{llnlElementalNavy}
\myth{Gen/Storage type} &
\myth{Plant lifetime (years)} &
\myth{Unit capacity (MW)} &
\myth{Annualized locational capital cost (2024 \$/MWy)} &
\myth{Locational FOM cost (2024 \$/MWy)} &
\myth{Locational VOM cost (2024 \$/MWh)} &
\myth{Source} \\
1x1 Combined Cycle (H-Class)           & 40 & 522.2            & 59107.38  & 15397.35 & 2.00 & DOE-FE\textsuperscript{1} \\
1x1 Combined Cycle (F-Class)         & 40 & 348.7            & 67215.43  & 19219.32 & 2.25 & DOE-FE\textsuperscript{1} \\
Large Simple Cycle (F-Class)         & 40 & 212.9            & 53535.52  & 14232.85 & 6.46 & DOE-FE\textsuperscript{1} \\
Small Simple Cycle (1 large aero)     & 40 & 103.9            & 81282.88  & 38418.73 & 3.74 & DOE-FE\textsuperscript{1} \\
Small Simple Cycle (2 small aero)     & 40 & 93.9             & 91892.67  & 42877.69 & 3.17 & DOE-FE\textsuperscript{1} \\
Small Simple Cycle (6 large RICE)     & 30 & 111.5            & 100970.29 & 37632.44 & 3.63 & DOE-FE\textsuperscript{1} \\
Small Simple Cycle (12 small RICE)    & 30 & 111.0            & 98857.86  & 37393.57 & 3.58 & DOE-FE\textsuperscript{1} \\
Onshore Wind\textsuperscript{2}           & 25 & 200.0            & 92400.88  & 35719.80 & 0.00 & AEO 2025 \cite{us_energy_information_administration_annual_2025} \\
Utility-scale Solar PV\textsuperscript{2} & 35 & 100.0            & 70590.04  & 24784.05 & 0.00 & AEO 2025 \cite{us_energy_information_administration_annual_2025} \\
BESS\textsuperscript{3}                   & 20 & 150.0, 4 hr      & 122022.17 & 41637.20 & 0.00 & AEO 2025 \cite{us_energy_information_administration_annual_2025} \\
\hline
\end{tabular}

\vspace{0.5em}
\raggedright
\textsuperscript{1} Conversations with DOE Office of Fossil Energy (DOE-FE).  Based on NETL report \cite{netl_report}.\\
\textsuperscript{2} Construction of fractional quantities is allowed, so unit capacity is only notional.\\
\textsuperscript{3} Battery energy storage system configured as 150 MW with 4 hours of storage.\\

\end{table*}

\subsection{Cost Assumptions for Fuel Costs}\label{sec:fuel_costs}
As discussed in Section \ref{sec:fuel_supply_data}, the location-specific fuel costs for San Juan, Palo Seco, Costa Sur, EcoElectrica, and the new Energiza project are relative to Henry Hub (HH) prices. We use the following HH prices for the existing (2024) and future (2030) cases:
\begin{center}
    \begin{tabular}{c|l}
  Existing & \$2.24/MMBtu\\
Future &\$3.15/MMBtu\\
 \end{tabular}
\end{center}
where the existing price is the historical HH average price from May-July 2024, and the future price is taken as the projected HH price in 2030 from EIA's most recent AEO, scaled by a seasonality factor of 1.023 for the May to July season based on the 2024 HH data.  Note that these are all in 2024 dollars.

Then we use the following equations (in \$/MMBtu), as provided by DOE, to determine the location-specific NG prices in Puerto Rico:
\begin{center}
    \begin{tabular}{>{\raggedright\arraybackslash}p{3.3cm} |l}
     SJ &6.5 + 1.15*HH\\
        PS (present - trucks)&7.60 + 1.15*HH + \$250/truck\\
      PS (future - pipe)& 7.60 + 1.15*HH\\
   EcoElectrica, CS& 5.5 + 1.15*HH\\
Energiza&7.95 + 1.15*HH
    \end{tabular}
\end{center}
where we consider the trucks to hold 860 MMBtu/truck, yielding an additional cost of \$0.29/MMBtu for trucked fuel.

We also obtained location-specific bunker fuel prices for San Juan, Palo Seco, Aguirre, and Costa Sur from the PREPA PI dashboard, from a snapshot of prices taken in October 2025.  Those are given as:
\begin{center}
    \begin{tabular}{c|l}
    San Juan & \$12.67/MMBtu\\
    Palo Seco& \$12.65/MMBtu\\
    Aguirre& \$13.63/MMBtu\\
    Costa Sur& \$12.98/MMBtu\\
 \end{tabular}
\end{center}

The diesel price of \$17/MMBtu was obtained from DOE and corresponds to December 2025 oil prices.  For the remaining fuel costs, including filling in NG prices for new units at Cambalache, Aguirre, etc., we used the LUMA December 2024 report \cite{luma_monthly_2025}.  Those are not location-specific and are given as:
\begin{center}
    \begin{tabular}{c|l}
     NG&\$11.90/MMBtu\\
   Bunker&\$14.30/MMBtu\\
  Coal& \$7.49/MMBtu\\
    \end{tabular}
\end{center}

\end{document}